\documentclass[11pt]{article}
\usepackage{graphicx} % Required for inserting images
\IfFileExists{bbm.sty}{\usepackage{bbm}}{}
\usepackage{bm}
\usepackage[margin=1in]{geometry}
\usepackage{color}
\usepackage{amsfonts}
\IfFileExists{algorithm.sty}{\usepackage{algorithm}}{}
\IfFileExists{algorithmic.sty}{\usepackage{algorithmic}}{}
\usepackage{latexsym, amssymb, amsmath, amscd, amsthm, amsxtra}
\usepackage{mathtools}
\usepackage{enumerate}
\usepackage[all]{xy}
\usepackage{mathrsfs}
\usepackage{fancyhdr}
\usepackage{listings}
\usepackage{hyperref}
\usepackage{enumitem}
\usepackage{multirow}
\usepackage{tikz}

\newcommand{\dist}{\mathsf{dist}}
\newcommand{\dep}{\mathsf{dep}}
\newcommand{\Ber}{\mathsf{Ber}}
\newcommand{\Bin}{\mathsf{Bin}}

\newcommand{\e}{\mathrm e}

\newcommand{\TV}{\mathsf{TV}}

\newcommand{\pu}{p_{\mathsf u}}
\newcommand{\ps}{p_{\mathsf s}}
\newcommand{\mix}{\mathsf{mix}}

\newcommand{\Pb}{\mathbb P}

\newcommand{\Eb}{\mathbb E}
\newcommand{\cG}{\mathcal{G}}

\newtheorem{thm}{Theorem}[section]

\newtheorem{prop}[thm]{Proposition}

\newtheorem{lemma}[thm]{Lemma}

\newtheorem{defn}[thm]{Definition}

\numberwithin{equation}{section}



\hypersetup{colorlinks=true,linkcolor=blue}

\title{Near linear mixing of random-cluster Glauber dynamics on random regular graphs}

\author{ Zhangsong Li\thanks{School of Mathematical Sciences, Peking University. Email: \textit{ramblerlzs@pku.edu.cn}} }

\date{\today}
\begin{document}
\maketitle

\begin{abstract}
We study the random-cluster model on locally treelike graphs at low temperatures. This is a model of dependent percolation parametrized by an edge probability $p \in (0,1)$ and a clustering weight $q \in [1,\infty)$, generalizing independent Bernoulli percolation ($q=1$) and closely related to the classical ferromagnetic Ising and Potts spin systems at integer $q$. For $q>2$, at parameter $p$ below the tree uniqueness threshold $p_{\mathsf u}(q,d)$, it is expected that sampling is easy and local Markov chains mix rapidly on all graphs with degree bounded by $d$. On typical graphs, the same is predicted at $p>p_{\mathsf s}(q,d)$, where $p_{\mathsf s}(q,d)$ is a second uniqueness transition point on the $d$-regular wired tree. Our result confirms this prediction for random regular graphs, showing that the random-cluster Glauber dynamics mix in time $O(n^{1+o(1)})$ on a random $d$-regular graph for all fixed $q>2$ and $p>p_{\mathsf s}$. 

For comparison, the best prior sampling results for these models that held for all $p>p_{\mathsf s}$ required $q\gtrsim \log d$ \cite{BGPZ26}. The key difficulty is that simply utilizing the contraction of cavity messages on $d$-regular trees does not suffice to cover the entire regime $q>2$. To address this issue, motivated by \cite{LS26}, we develop a careful bootstrap argument to show the shrinkage of disagreement under grand coupling. The core ideas of the proof were generated using GPT-6 and the Pharos system; the proof was then refined and simplified by the author.
\end{abstract}

\tableofcontents

\section{Introduction}{\label{sec:intro}}

We study the random-cluster model and its local Glauber dynamics Markov chain on locally treelike graphs in the low-temperature regime. The random-cluster model is an extensively studied model of random subgraphs both of interest in its own right and for its deep connections to the ferromagnetic Ising and Potts models \cite{Grimmett06}. Formally, on a finite graph $G= (V,E)$, the random-cluster model is the probability distribution over $A \subseteq E$ given by
\begin{equation}{\label{eq:def-random-cluster}}
    \mu_{G,p,q}(A) \propto p^{|A|} (1-p)^{|E \setminus A|} q^{c(A)},
\end{equation}
where $c(A)=c_G(A)$ counts the number of connected components of the subgraph $(V,A)$, and $p\in (0,1), q\in [1,\infty)$ are model parameters. At $q=1$, the distribution \eqref{eq:def-random-cluster} is simply the independent edge percolation, while at integer $q \ge 2$ it admits a coupling to the classical ferromagnetic Ising and Potts models. This coupling allows samples generated for one model to be converted into samples from the other and is also at the heart of the famous Swendsen-Wang Markov chain, one of the most widely used algorithms for sampling from these distributions. We study the random-cluster model when $q>1$ is any real number. The $q<1$ regime is also of interest, but qualitatively very different because the positive-association and monotonicity arguments used for $q\ge1$ no longer apply; see \cite{ALGV24} for important work on sampling from the $q<1$ random-cluster model at all $p$ on all graphs.

The random-cluster model generically undergoes a phase transition as the parameter $p$ is varied. At low values of $p$ (high temperature), correlations between edges decay exponentially fast in their distance and no large connected components appear. In contrast, at large values of $p$ (low temperature), a unique ``giant'' component of linear size emerges and all other cluster sizes have exponential tails on their diameters. For $q>2$, however, there can be a \emph{coexistence} regime of these two different phases at intermediate values of $p$. The specifics of this phase transition depend on the underlying graph. On the integer lattice $\mathbb Z^d$, the subcritical and supercritical regimes are separated by a single critical point $p=p_{\mathsf c}$ \cite{DCRT19}. On the wired regular tree, and in the corresponding metastability picture for random regular graphs, the coexistence regime can instead occupy an entire interval of parameter values, denoted by $(\pu,\ps)$. This poses special algorithmic and mathematical challenges as the transition points $\pu$ and $\ps$ do not align with the point at which a typical sample from the distribution transitions from high vs. low temperature behavior. Notably, the coexistence regime $(\pu,\ps)$ is characterized by the slow convergence of natural Markov chains and underlies the computational hardness of approximate counting and sampling for the random-cluster and Potts models for $q>2$ \cite{GJ12, GSVY16}. 

Focusing on Markov chain sampling, we recall the definition of the heat-bath Glauber dynamics for the random-cluster model, which we henceforth refer to as \emph{random-cluster dynamics}. In discrete-time, this is the Markov chain which from configuration $X_t \subseteq E$, generates $X_{t+1}$ by
\begin{itemize}
    \item Picking an edge $e \in E$ uniformly at random;
    \item Setting $X_{t+1}$ to be
    \begin{equation}{\label{eq:transition-random-cluster-dynamics}}
        X_{t+1} = X_t \cup \{ e \} \mbox{ with probability }
        \begin{cases}
            \widehat{p} := \frac{p}{q(1-p)+p}, & e \mbox{ is a ``cut edge'' in } (V,X_t\cup\{e\}); \\
            p, & \mbox{otherwise};
        \end{cases}
    \end{equation}
    and $X_{t+1}=X_t\setminus\{e\}$ otherwise.
\end{itemize}
The mixing time of a Markov chain is the number of steps it takes it to be close in total-variation to the stationary distribution $\mu_{G,p,q}$ from a worst-case initialization; namely,
\begin{equation}{\label{eq:def-mixing-time}}
    T_{\mix}(G,p,q) = \min\left\{ t \ge 0: \max_{X_0} \TV\big( \Pb(X_t \in \cdot \mid X_0),\mu_{G,p,q} \big) \le \frac{1}{4} \right\}.
\end{equation}
When $q \in (1,2]$, the mixing time of the random-cluster dynamics is expected to always be polynomial on any input graph for any $p \in (0,1)$; this was shown at $q=2$ in \cite{GJ18}. When $q>2$, the expectation is that at high temperatures $p<\pu$ the random-cluster dynamics mix quickly, slow down exponentially in the intermediate coexistence regime $p \in (\pu,\ps)$, but interestingly become rapidly mixing again as soon as $p>\ps$. This trichotomy of behaviors has been fully characterized on the complete graph \cite{BGJ96, BS15, GSV15, GLP19}.

Beyond the complete graph, considerable work has focused on establishing this trichotomy for other non-amenable graph families, especially locally tree-like random graphs. Physicists have long conjectured that graphical models on these locally tree-like ensembles exhibit thermodynamic behavior very similar to those on regular lattices, especially with respect to phase transitions and the structure of Gibbs measures (e.g. \cite{ACTA73, CLR79, Thouless86}). At the same time, the absence of short loops and the recursive nature of their local neighborhoods often make these models more analytically tractable--and in many cases exactly solvable. Building on physics predictions, over the past fifteen years, rigorous studies have examined graphical models on locally tree-like graphs, focusing on various aspects including the limiting free-energy density \cite{DM10, DMS13, SS14, DGvdH14, DMSS14, Can19, BBC23, GGvdH+25, CvdH25} and the local structure of Gibbs measures \cite{MMS12, BD17, HJP23, BDS23, YP24, DZ25}.
However, the model behavior for general $q>1$ and $p \in (0,1)$ poses significant mathematical challenges and remains not fully understood at low temperatures, neither in terms of its decorrelation properties nor its algorithmic tractability. 

On the $d$-regular tree, \cite{Haggstrom96} showed that for $p<\pu$, there is decay of correlations between the edges at the top of the tree and the boundary, and notably that there is a unique Gibbs random-cluster measure on the infinite $d$-regular wired tree. The follow-up work \cite{BG21} proved a mixing time of $\Theta(n\log n)$ on random $d$-regular graphs for every $q\ge1$ and $p<\pu$. Non-uniqueness and lack of correlation decay were also shown for the $p \in (\pu,\ps)$ regime on the wired tree in \cite{Haggstrom96}, where
\begin{equation}{\label{eq-def-ps}}
    \ps = \ps(q,d) = \frac{q}{d+q-2}.
\end{equation}
(The threshold $\pu$ has also been identified, but only implicitly as it does not admit a closed form.)
On random regular graphs, slow mixing when $p \in (\pu,\ps)$ was shown in \cite{COGG+23} for integer $q$, and in \cite{BBC23} for general real $q$. Thus, the remaining question is whether rapid mixing holds when $p>\ps$. In the low-temperature regime, \cite{BGPZ26} proved near-linear mixing throughout $p>\ps$ under the assumption $q\gtrsim\log d$. Our result removes the undesirable assumption $q \gtrsim \log d$, showing that the random-cluster dynamics on random $d$-regular graphs indeed mix in time $O(n^{1+o(1)})$ throughout $p>\ps$.
\begin{thm}{\label{thm:main}}
    Fix an integer $d \ge 3$, a real number $q>2$, and $\ps<p<1$. For every constant $\eta>0$, there exists a constant $C=C(p,q,d,\eta)$ such that 
    \[
        \Pb_{G \sim \cG_d(n)}\Big( T_{\mix}(G,p,q) \le Cn^{1+\eta} \Big) \to 1 \mbox{ as } n \to \infty.
    \]
    Here $\cG_d(n)$ is the uniform distribution over all simple $d$-regular graphs on $V=\{1,\ldots,n\}$.
\end{thm}

\subsection{Notation}{\label{subsec:notation}}

We record in this subsection some notation conventions. For any connected graph $H$, denote the tree excess of $H$ by $\tau(H)=|E(H)|-|V(H)|+1$. Given $u \in V(H)$, denote by $\mathcal{N}_H(u)$ the set of neighbors of $u$ in $H$, and denote by $\deg_H(u)$ the degree of $u$ in $H$. For two vertices $u,v\in V(H)$, we define $\dist_{H}(u,v)$ to be the graph distance. Denote by $\mathsf L(H)$ the set of leaves of $H$, and denote by $\mathcal C(H)$ the set of connected components of $H$. For a graph $H=(V,E)$ and a subset $U \subseteq V$, define $H[U] =(U,E_H(U))$ to be the vertex-induced subgraph of $H$ in $U$, where $E_H(U):=\{\{u,v\}\in E:u,v\in U\}$. Also define $\partial_H(U):=\{\{u,v\}\in E:u\in U,\ v\notin U\}$. For $v \in V(G)$ and an integer $r \ge 0$, set
\begin{equation}{\label{eq:def-local-ball}}
    B_{v,r}=\{ w: \dist_G(v,w) \le r \}, \qquad S_{v,r}=\{ w: \dist_G(v,w) = r \}, \qquad E_{v,r}= E(G[B_{v,r}]).
\end{equation}
For $A \subseteq E$, write 
\begin{equation}{\label{eq:def-fragment-cluster}}
    F_{v,r}(A) = (V,A \setminus E_{v,r}), \qquad C_{v,r}(A) = C_1(F_{v,r}(A)).
\end{equation}
Here $C_1(H)$ is the vertex set of the largest connected component of $H$ (size ties are broken by the least vertex label). It is clear that $F_{v,r}(A)$ induces a partition of $S_{v,r}$, where vertices belong to the same group if and only if they are connected in $F_{v,r}(A)$. Denote by $\zeta_{v,r}(A)$ this partition.

Define $\rho(\mathbb T)$ to be the root of a tree $\mathbb T$. Denote by $E_\rho(\mathbb T)$ the set of edges of $\mathbb T$ incident to the root. For a rooted tree $(T,o)$ and $u\in V(T)$, we define $\dep_{T}(u)$ to be the depth of $u$ in $T$, and define $\dep(T)$ to be the depth of $T$. For $u,v \in V(T)$, denote by $u \rightarrow v$ (or equivalently, $v \leftarrow u$) if $v$ is a child of $u$, and denote by $u \hookrightarrow v$ (or equivalently, $v \hookleftarrow u$) if $v$ is a descendant of $u$.

Throughout the remainder of the paper, we will fix $d \ge 3$, $q>2$, $\ps<p<1$, and $\eta>0$. All constant dependencies on $d,q,p,\eta$ will be suppressed. For any two positive sequences $\{a_n\}$ and $\{b_n\}$, we write interchangeably $a_n=O(b_n)$, $b_n=\Omega(a_n)$, $a_n\lesssim b_n$ and $b_n\gtrsim a_n$ if there exists a positive constant $c$ depending only on the fixed parameters such that $a_n/b_n\leq c$ holds for all $n$ (similarly we use the notation $O_h(\cdot)$ to allow additional dependence on $h$). We write $a_n=o(b_n)$, $b_n=\omega(a_n)$, $a_n\ll b_n$, and $b_n\gg a_n$ if $a_n/b_n\to 0$ as $n\to\infty$. We say a family $\cG$ of $d$-regular graphs is a typical family, if $\Pb_{G \sim \cG_d(n)}(G \in\cG)=1-o(1)$.

\section{Heuristics and proof overview}{\label{subsec:intuition}}

We explain the main ideas behind Theorem~\ref{thm:main}, emphasizing how our argument removes the restriction $q\gtrsim\log d$ in \cite{BGPZ26}. 

\subsection{A bootstrap contraction argument}{\label{subsec:overview-contraction}}

We begin by recalling the proof of \cite{BGPZ26}, explaining why it is limited to $q\gtrsim \log d$. The key step in the proof in \cite{BGPZ26} is to establish uniqueness and decay of correlations on the wired $d$-regular tree. At a high level, the reason $\ps$ from \eqref{eq-def-ps} is predicted as the threshold for uniqueness, is that when $p>\ps$, one has that $\widehat{p}$ defined in \eqref{eq:transition-random-cluster-dynamics} is bigger than $\frac{1}{d-1}$, and the random-cluster measure stochastically dominates a supercritical Bernoulli percolation on the tree. Thus, any random-cluster measure on the tree has a positive density of infinite components that are all connected together in the wired tree, and at $p>\ps$ any such positive density should be enough to make the model behave like it does when all of its leaves are wired together. To justify this intuition, \cite{BGPZ26} exploit the following approach. First, for each $e \in E$, they choose a local, tree-like neighborhood $H_e$. They then compare, on each local neighborhood $H_e$, a sufficiently wired boundary condition obtained after a burn-in period with the fully wired boundary condition. The cavity recursion controls the difference between their root messages, and hence between their root-edge marginals. To make the total coupling error vanish by a union bound over the $O(n)$ edges, these marginal differences must be $o(n^{-1})$. Since the available neighborhoods have depth $O(\log n/\log(d-1))$, this requires the cavity contraction factor to be a sufficiently small inverse power of $d-1$. Obtaining such a rate throughout $p>\ps$, including parameters arbitrarily close to $\ps$, forces $q\gtrsim\log d$ within this approach; see \cite[\S4--5]{BGPZ26}.

To overcome this limitation, our first new ingredient is a bootstrap argument for the disagreement size. We show that, after burn-in, the expected number of disagreements decreases by a constant factor over each time window, up to a vanishing error. We first formulate the desired contraction under the grand coupling.
%This requires controlling the cumulative effect of the boundary differences on all local trees by the number of existing disagreements. The local estimate must therefore retain the dependence on the actual boundary changes, which is lost when both boundaries are compared separately with the fully wired condition. 
\begin{defn}[Grand coupling]\label{def:grand-coupling}
    Let $(e_t)_{t\ge1}$ be independent uniform edges of $E$, and let $(U_t)_{t\ge1}$ be independent $\mathrm{Unif}[0,1]$ random variables, independent of the update edges. For $B\subseteq E$ and $e\in E$, put
    \[
        \theta_e(B)=
        \begin{cases}
            p, & \text{if the endpoints of $e$ are connected in }(V,B\setminus\{e\}),\\
            \widehat p, & \text{otherwise}.
        \end{cases}
    \]
    For each initial configuration $A\subseteq E$, set $X_0^A=A$ and
    \[
        X_t^A=
        \begin{cases}
            X_{t-1}^A\cup\{e_t\}, & U_t\le\theta_{e_t}(X_{t-1}^A),\\
            X_{t-1}^A\setminus\{e_t\}, & U_t>\theta_{e_t}(X_{t-1}^A).
        \end{cases}
    \]
    Using the same update edges and uniform variables for all $A$ gives the grand coupling. Since $q\ge1$, this coupling is monotone: if $A\subseteq A'$, then $X_t^A\subseteq X_t^{A'}$ for all $t\ge0$.
\end{defn}

Let $X_t^-$ and $X_t^+$ be the chains started from $\emptyset$ and $E$, respectively. In addition, draw $X_0\sim\mu_{G,p,q}$ independently of the update sequence, and let $X_t$ be the corresponding stationary chain under the grand coupling. Then
\[
    X_t^-\subseteq X_t\subseteq X_t^+\qquad(t\ge0).
\]
Define the expected number of disagreements by
\begin{equation}\label{eq:def-disagreement-size}
    D_t:=\mathbb E\bigl[|X_t^+\setminus X_t^-|\bigr].
\end{equation}
The coupling inequality gives
\[
    \max_{A\subseteq E}\Big\{ \TV\bigl(\Pb(X_t^A\in\cdot),\mu_{G,p,q}\bigr) \Big\} \le\Pb(X_t^-\ne X_t^+)\le D_t.
\]
Our goal is to find deterministic times $T_0=O(n\log n)$ and $T=O(n^{1+\eta/2})$ such that, for some $c>0$,
\begin{equation}\label{eq:overview-contraction}
    D_{t+T}\le\frac14D_t+n^{-c} \qquad\text{for every deterministic integer }t\ge T_0.
\end{equation}
Clearly, once \eqref{eq:overview-contraction} is settled, we can derive Theorem~\ref{thm:main} by iterating \eqref{eq:overview-contraction} up to $\log n$ steps; see \S\ref{subsec:overview-finishing} for further details.

\subsection{Local neighborhood comparison}{\label{subsec:overview-local-comparison}}

We now explain how to derive \eqref{eq:overview-contraction}. Suppose at time $t$ we have $X^{-}_t =A^{-}, X_t^{+}=A^{+}$ for some $A^{-} \subseteq A^{+} \subseteq E$. For $A \subseteq E$, let $\xi_e(A)$ be the partition of $S_e$ induced by connectivity in $(V,A\setminus E(H_e))$. For each local neighborhood $H_e$, the two different configurations may lead to two different boundary conditions $\xi_{e}(A^{\pm})$. A natural attempt is to control the cumulative effect of the boundary differences on all local trees by the number of existing disagreements. The local estimate must therefore retain the dependence on the actual boundary changes, which is lost when both boundaries are compared separately with the fully wired condition. We need two ingredients to carry out such an argument:
\begin{itemize}
    \item A comparison result which controls the difference of root messages under the boundary condition $\xi_{e}(A^{+})$ and $\xi_{e}(A^{-})$.
    \item A local mixing argument which shows that, restricted to the local neighborhood $H_e$, the random-cluster dynamics near the root converge rapidly to the equilibrium.
\end{itemize}
However, there are several obstacles to carrying out these two arguments. The first difficulty is that rigorously speaking, the local neighborhood $H_e$ is only tree-like instead of always being precisely a tree, which makes the precise cavity calculations more complicated or even untractable. We will overcome this by the following more delicate choice of the neighborhood $H_e$. Choose an endpoint $o_e$ of each edge $e$ by a deterministic rule depending only on $G$. For constants $0<\alpha<\beta<1/4$ to be chosen later in \eqref{eq:choice-alpha-beta}, set
\begin{equation}\label{eq:def-radii-r-R}
    r=\lfloor\alpha\log_{d-1} n\rfloor, \qquad R=\lfloor\beta\log_{d-1} n\rfloor.
\end{equation}
Write $B_{v,h}=\{w:\dist_G(v,w)\le h\}$ and $S_{v,h}=\{w:\dist_G(v,w)=h\}$. Define the exceptional edges by
\begin{equation}\label{eq:def-bad-edges-mathcal-B}
    \mathcal B=\{e\in E:G[B_{o_e,r}]\text{ is not a tree}\},
\end{equation}
and use the comparison balls
\begin{equation}\label{eq:def-local-ball-H-e}
    H_e=
    \begin{cases}
        G[B_{o_e,r}], & e\notin\mathcal B,\\
        G[B_{o_e,R}], & e\in\mathcal B.
    \end{cases}
    \qquad 
    S_e=
    \begin{cases}
        S_{o_e,r}, & e\notin\mathcal B,\\
        S_{o_e,R}, & e\in\mathcal B.
    \end{cases}
\end{equation}
Intuitively, we will first choose the $r$-neighborhood of the edge $e$, and we will enlarge the depth of the neighborhood if its original neighborhood is not a tree. Since a typical random regular graph is tangle-free (see \cite{Friedman08}), the enlarged part of the neighborhood does not create cycles, and it yields a sharper contraction of cavity messages which enables us to beat the enumeration of $\mathcal B$. We now introduce the following notation, which reduces the problem to studying the behavior of dynamics restricted to each local neighborhood.
\begin{defn}[Local chains]\label{def:local-chain}
    Given $A^{-} \subseteq A^{+}\subseteq E$ and $e\in E$, start $\widetilde X_0^{e,-}=A^{-}$ and $\widetilde X_0^{e,+}=A^{+}$. Use the grand-coupling updates on $G$, retaining an update if its selected edge belongs to $E(H_e)$ and otherwise doing nothing. Thus the two exteriors remain fixed at $A^{-} \setminus E(H_e)$ and $A^{+} \setminus E(H_e)$. Note that the stationary laws of the two local chains, restricted to $E(H_e)$, are therefore 
    \[
        \mu_{H_e,p,q}^{\xi_e(A^-)} \quad \text{and} \quad \mu_{H_e,p,q}^{\xi_e(A^+)}.
    \]
\end{defn}
For a common time window $T \ge 0$, define
\begin{equation}\label{eq:q-e-X-pm}
    q_e(A^-,A^+):=\Pb\Big( e \in \widetilde X_T^{e,+} \setminus \widetilde X_T^{e,-} \Big) =\Pb\Big( e \in \widetilde X_T^{e,+} \Big) - \Pb\Big( e \in \widetilde X_T^{e,-} \Big).
\end{equation}
We then need to bound the cumulative effect of $q_e(A^-,A^+)$ by the number of existing disagreements $|A^+ \setminus A^-|$. 
However, two more difficulties arise here. First, when $p>\ps$ we do not have a good mixing estimate for random-cluster dynamics on trees with a general boundary condition, even under the assumption that the boundary contains a ``sufficiently wired'' block, so in general we cannot assume that $\widetilde X_t^{e,+}$, $\widetilde X_t^{e,-}$ converge to their equilibrium rapidly. Second, even if we know $\widetilde X_T^{e,+}$, $\widetilde X_T^{e,-}$ are close to their stationary laws, comparing the difference of these two equilibria near the root is still challenging for general boundary conditions.

To handle these two issues, our observation is that we do have a good mixing estimate and a good comparison if the two boundary conditions $\xi_e(A^-), \xi_e(A^+)$ satisfy certain desirable conditions. Roughly speaking, the local mixing and comparison can be proved if these two boundary conditions have a bounded number of non-singleton blocks, and one of the blocks is ``sufficiently dense''; see Lemmas~\ref{lem:multi-block-boundary-comparison} and \ref{lem:spec-gap-small-excess}. Our second innovation is to introduce \emph{comparison boundary conditions}
\[
    \zeta_e^- \preceq \xi_e(A^-) \preceq \xi_e(A^+) \preceq \zeta_e^+
\]
such that $\zeta_e^{\pm}$ satisfy these desirable properties. Thus, by the monotonicity of random-cluster measure, at a common time $T=O(n^{1+\eta/2})$, we are able to show that
\[
    q_e(A^-,A^+)\le\Delta_e+O(n^{-12}),\qquad \Delta_e:=\mu_{H_e,p,q}^{\zeta_e^+}(e\mbox{ open}) -\mu_{H_e,p,q}^{\zeta_e^-}(e\mbox{ open}).
\]
Finally, on tree neighborhoods, changes in the comparison boundaries can be charged to the disagreement set $A^+\setminus A^-$. On the few exceptional neighborhoods, a larger radius gives a small absolute error. Combining these estimates yields
\[
    \sum_{e\in E}q_e(A^-,A^+)\le\frac14|A^+\setminus A^-|+n^{-c_0}
\]
whenever the two configurations contain a suitable product configuration and enclose a suitable stationary configuration. A standard censoring argument then enables us to bound $D_{t+T}$ by $D_t$.

These estimates must hold uniformly over all pairs $A^-\subseteq A\subseteq A^+$ containing a suitable product configuration $Y$, since the configurations generated by the dynamics are dependent on the reference configuration $A$. We now explain how the grand coupling provides these reference configurations. Let $Y_0\sim\Ber(\widehat p)^{\otimes E}$, independently of $X_0$ and the update sequence. At time $t$, refresh the selected edge $e_t$ of $Y$ to $\mathbf1_{\{U_t\le\widehat p\}}$, leaving all other edges unchanged. For every deterministic $t$,
\[
    Y_t\sim\Ber(\widehat p)^{\otimes E},\qquad X_t\sim\mu_{G,p,q}.
\]
Set
\begin{equation}\label{eq:def-E-1}
    T_0=\lceil100dn\log n\rceil, \qquad \mathcal E_{\diamond}=\{E\subseteq\{e_1,\ldots,e_{T_0}\}\}.
\end{equation}
Coupon collection gives $\Pb(\mathcal E_{\diamond}^c)\le n^{-20}$ for all sufficiently large $n$. Since every random-cluster opening probability is at least $\widehat p$, on $\mathcal E_{\diamond}$ we have
\[
    Y_t\subseteq X_t^-\subseteq X_t\subseteq X_t^+
    \qquad\text{for all }t\ge T_0.
\]
The following property expresses the uniform estimate needed from a pair of comparison configurations.
\begin{defn}[Admissible pair]\label{def:admissible-realization}
    Fix the balls in \eqref{eq:def-local-ball-H-e}, the time window $T$, and a constant $c_0>0$. A pair $Y\subseteq A\subseteq E$ is admissible if, for every $(A^-,A^+)$ satisfying $Y\subseteq A^{-}\subseteq A\subseteq A^+\subseteq E$,
    \[
        \sum_{e\in E}q_e(A^-,A^+)\le\frac14|A^+\setminus A^-|+n^{-c_0}.
    \]
\end{defn}
The main step in our proof is to show the following result.
\begin{prop}{\label{prop:admissible}}
    There are constants $c_0,C>0$ and a typical family $\cG_{\ref{prop:admissible}}$, such that for every $G\in\mathcal G_{\ref{prop:admissible}}$ one can choose a deterministic integer $T\le Cn^{1+\eta/2}$ for which, at every deterministic integer $t\ge T_0$,
    \begin{equation}\label{eq:overview-good-pair}
        \Pb_G\Big( \mathcal E_\diamond\cap\{(Y_t,X_t)\text{ is admissible}\} \Big) \ge 1-Cn^{-20}.
    \end{equation} 
\end{prop}

\subsection{Finishing the proof}{\label{subsec:overview-finishing}}

We now deduce Theorem~\ref{thm:main} from Proposition~\ref{prop:admissible}.
\begin{proof}[Proof of Theorem~\ref{thm:main} assuming Proposition~\ref{prop:admissible}]
    We first compare the original dynamics with the local chains. Fix a deterministic integer $t\ge T_0$ and an edge $e\in E$. Starting from $\emptyset$ and $E$ at time zero, define censored chains $\widehat X_s^{e,-}$ and $\widehat X_s^{e,+}$ by retaining every original update through time $t$, and retaining only updates in $E(H_e)$ during the interval $(t,t+T]$. In particular, $\widehat X_t^{e,\pm}=X_t^\pm$. Conditional on the complete sequence of selected edges through time $t+T$, the censored schedule is a fixed subsequence of the original schedule. The Peres--Winkler censoring inequality \cite[Theorem~1.1]{PW13}, applied from the two extremal initial configurations, therefore gives
\[
    \Pb\Big( e\in X_{t+T}^+ \Big) \le \Pb\Big( e\in\widehat X_{t+T}^{e,+} \Big), \qquad \Pb\Big( e\in X_{t+T}^- \Big) \ge \Pb\Big(e\in\widehat X_{t+T}^{e,-} \Big).
\]
Conditional on the full history of the coupled processes through time $t$, the next $T$ steps of these censored chains are exactly the local chains in Definition~\ref{def:local-chain}, started from $X_t^-$ and $X_t^+$. The future update edges and uniforms are independent of this history. Hence, by \eqref{eq:q-e-X-pm},
\[
    \Pb\Big( e\in X_{t+T}^+ \Big) - \Pb\Big( e\in X_{t+T}^- \Big) \le\mathbb E\Big[ q_e(X_t^-,X_t^+) \Big].
\]
Summing over $e$ and using monotonicity yields
\[
    D_{t+T}\le\mathbb E\left[ \sum_{e\in E}q_e(X_t^-,X_t^+) \right].
\]
Let $\mathcal A_t:=\mathcal E_\diamond\cap\{(Y_t,X_t)\text{ is admissible}\}$.
On $\mathcal A_t$, the inclusion $Y_t\subseteq X_t^-\subseteq X_t\subseteq X_t^+$ allows us to apply Definition~\ref{def:admissible-realization} with $A^-=X_t^-$ and $A^+=X_t^+$. On its complement, we use $0\le q_e\le1$. Proposition~\ref{prop:admissible} thus gives
\[
\begin{aligned}
    D_{t+T}
    &\le\frac14\mathbb E\bigl[\mathbf1_{\mathcal A_t}|X_t^+\setminus X_t^-|\bigr]
           +n^{-c_0}\Pb(\mathcal A_t)+dn\Pb(\mathcal A_t^c)\\
    &\le\frac14D_t+n^{-c_0}+C_1n^{-19}
     =\frac14D_t+\delta_n,
    \qquad \delta_n:=n^{-c_0}+C_1n^{-19}.
\end{aligned}
\]
For $c=\tfrac12\min\{c_0,19\}$ and all sufficiently large $n$, we have $\delta_n\le n^{-c}$, proving \eqref{eq:overview-contraction}. 

Finally, set $j_n=\lceil2\log_4 n\rceil$. Iterating the recurrence at $T_0,T_0+T,\ldots,T_0+(j_n-1)T$, and using $D_{T_0}\le |E(G)|$, we obtain
\[
    D_{T_0+j_nT}
    \le4^{-j_n}|E(G)|+\delta_n\sum_{i=0}^{j_n-1}4^{-i}
    \le\frac{d}{2n}+\frac43\delta_n=o(1).
\]
For every initial configuration $A\subseteq E$, the grand coupling places both $X_u^A$ and the stationary chain $X_u$ between $X_u^-$ and $X_u^+$. Therefore, at the deterministic time $u_n=T_0+j_nT$,
\[
    \max_{A\subseteq E}\Big\{ \TV\bigl(\Pb(X_{u_n}^A\in\cdot),\mu_{G,p,q}\bigr) \Big\}
    \le\Pb(X_{u_n}^-\ne X_{u_n}^+)\le D_{u_n}\le\frac14
\]
for all sufficiently large $n$. It follows that $T_{\mix}\le u_n$. The bounds on $T_0$ and $T$ give
\[
    u_n\le C'\bigl(n\log n+n^{1+\eta/2}\log n\bigr)
    \le C''n^{1+\eta}. \qedhere
\]
\end{proof}

\section{Structure of random regular graphs}{\label{sec:random-regular-graph-structure}}

We now formalize the desirable expansion properties we require for typical $G \sim\mathcal G_d(n)$, as described in Section~\ref{subsec:intuition}.
\begin{defn}{\label{def:typical-graph}}
    Put 
    \begin{equation}\label{eq:structural-parameters}
        \epsilon=\min\left\{ \frac{(d-1)\widehat{p}-1}{100},1 \right\}, \qquad \delta=(100d)^{-2/\epsilon}.
    \end{equation}
    Define $\mathcal G_{\ref{def:typical-graph}}$ to be the set of connected, $d$-regular graphs $G=(V,E)$ on $n$ vertices such that the following conditions hold.
    \begin{itemize}
        \item[(1)] For all $v \in V$, we have $\tau(G[B_{v,R}]) \le 1$, where $R$ is defined in \eqref{eq:def-radii-r-R}. 
        \item[(2)] We have $\#\{ v \in V: \tau(G[B_{v,r}]) \ge 1 \} \le n^{3\alpha}$.
        \item[(3)] For all $U \subseteq V$ such that $1 \le |U| \le \frac{n}{2}$, we have $|\partial_G(U)| \ge 0.01d|U|$.
        \item[(4)] For all $U \subseteq V$ with $1 \le |U| \le \delta n$, we have $|E(G[U])| \le (1+\epsilon)|U|$.
        \item[(5)] Recall \eqref{eq:def-fragment-cluster}. For all $v \in V$ and for all $H \subseteq F_{v,r}(E)$ such that $H$ is connected and $|V(H)| \ge \log n$, we have
        \[
            \big| V(H) \cap S_{v,r} \big| \le \frac{\epsilon}{d}|V(H)|, \qquad \sum_{x \in V(H)} \Big( d-\mathsf{deg}_{F_{v,r}(E)}(x) \Big) \le \epsilon |V(H)|.
        \]
    \end{itemize}
\end{defn}
The main goal of this section is to show the following structural results on a random regular graph $G \sim \mathcal G_d(n)$.
\begin{prop}{\label{prop:typical-graph}}
    Suppose $G \sim\mathcal G_d(n)$ where $d \ge 3$, then $\Pb(G \in \mathcal G_{\ref{def:typical-graph}})=1-o(1)$.
\end{prop}
The rest of this section is devoted to the proof of Proposition~\ref{prop:typical-graph}. By \cite[Fact~2.17]{GSS25}, Item (1) holds for $G \sim\mathcal G_d(n)$ with probability $1-o(1)$. We then prove that Item~(2) holds for $G \sim\mathcal G_d(n)$ with probability $1-o(1)$ via the following lemma.
\begin{lemma}{\label{lem:num-small-cycles}}
    Fix an integer $d\ge3$ and constants $0<\alpha<\beta<1/4$, and recall $r,R$ from \eqref{eq:def-radii-r-R}. If $G\sim\mathcal G_d(n)$, then, with probability $1-o(1)$,
    \[
        \#\Big\{ v\in V:\tau(G[B_{v,r}])\ge1 \Big\}\le n^{3\alpha}.
    \]
\end{lemma}
\begin{proof}
    It suffices to prove the assertion in the configuration model. Fix $v\in V$ and explore the ball around $v$ by breadth-first search, pairing all half-edges based at vertices of distance less than $r$. Once $B_{v,r}$ has been determined, also pair every remaining half-edge based at a vertex of $S_{v,r}$. This last step reveals all edges within the outer sphere, so the exploration determines the entire induced subgraph $G[B_{v,r}]$. The degree bound gives
    \[
        |B_{v,r}|\le1+d\sum_{i=0}^{r-1} (d-1)^i\le C_d (d-1)^r.
    \]
    Consequently, the exploration makes at most $M_r:=\lceil C_d (d-1)^r \rceil$ pairing steps and touches at most $M_r+1$ vertices, including exterior vertices reached when pairing the sphere half-edges. Call a pairing a collision if its partner belongs to a vertex already touched by the exploration. Conditional on the preceding exploration, the partner is uniform among the remaining unmatched half-edges. Since $(d-1)^r \le n^\alpha=o(n)$, its conditional probability of being a collision is at most
    \[
        \frac{d(M_r+1)}{dn-2M_r-1} \le \frac{C_d'(d-1)^r}{n}
    \]
    for all sufficiently large $n$. A union bound therefore gives
    \begin{equation}\label{eq:ball-collision-bound}
        \Pb\Big( \tau(G[B_{v,r}])\ge 1 \Big) \le M_r \cdot \frac{C_d'(d-1)^r}{n} \le \frac{ C_{d}(d-1)^{2r}}{n} .
    \end{equation}
    Let $N_r=\#\{v\in V:\tau(G[B_{v,r}])\ge1\}$. Applying \eqref{eq:ball-collision-bound} and summing over $v$ gives
    \[
        \Eb[N_r]\le C_{d} (d-1)^{2r}\le C_{d} n^{2\alpha}.
    \]
    Markov's inequality implies
    \[
        \Pb(N_r>n^{3\alpha})\le C_{d}n^{-\alpha}=o(1).
    \]
    Conditioning on simplicity proves the assertion for $G\sim\mathcal G_d(n)$.
\end{proof}
We then prove that Item~(3) holds for $G \sim\mathcal G_d(n)$ with probability $1-o(1)$ via the following lemma.
\begin{lemma}{\label{lem:edge-expansion}}
    Fix an integer $d\ge3$, and let $G\sim\mathcal G_d(n)$. Then, with probability $1-o(1)$, we have 
    \[
        |\partial_G(U)|\ge 0.01d|U| \quad \text{for every} \quad U \subseteq V, \ 1\le|U|\le\frac{n}{2}.
    \]
\end{lemma}
\begin{proof}
    Let $M_G$ be the adjacency matrix of $G$, with eigenvalues $d=\lambda_1\ge\lambda_2\ge\cdots\ge\lambda_n$. Friedman's second-eigenvalue theorem (see, e.g., \cite{Friedman08} or \cite[Theorem~1]{Bor20}) therefore implies that, with probability $1-o(1)$,
    \[
        \lambda_2 \le2\sqrt{d-1}+0.01d < d-0.02d.
    \]
    Fix a graph in this event and a set $U\subseteq V$ with $u:=|U|\in[1,n/2]$. Put $f=\mathbf1_U-(u/n)\mathbf1$, where $\mathbf1_U$ is the indicator vector of $U$ and $\mathbf1$ is the all-ones vector. Then $f\perp\mathbf1$ and $\|f\|_2^2=u(1-\frac un)$.
    Since $G$ is $d$-regular, its Laplacian is $dI-M_G$, and
    \[
        f^{\top}(dI-M_G)f =\sum_{\{x,y\}\in E}(f(x)-f(y))^2 =|\partial_G(U)|.
    \]
    The variational characterization of $\lambda_2$ therefore gives
    \[
        |\partial_G(U)| \ge(d-\lambda_2)\|f\|_2^2 \ge 0.02du\left(1-\frac un\right) \ge 0.01du. \qedhere
    \]
\end{proof}
We then prove that Item~(4) holds for $G \sim\mathcal G_d(n)$ with probability $1-o(1)$, as incorporated in the following lemma.
\begin{lemma}{\label{lem:subgraph-density}}
    Recall \eqref{eq:structural-parameters}. If $G\sim\mathcal G_d(n)$, then, with probability $1-o(1)$, we have
    \[
        |E(G[U])| \le (1+\epsilon)|U| \quad \text{for every} \quad U \subseteq V, \ 1 \le |U| \le \delta n.
    \]
\end{lemma}
\begin{proof}
    It suffices to work with the configuration model. Recall that $0<\epsilon\le1$ and $\delta\le(100d)^{-2}<1/2$. Fix a set $U$ with $|U|=t\le\delta n$, and put $j=\lfloor(1+\epsilon)t\rfloor+1$. If $2j>dt$, the event under consideration is impossible, so assume $2j\le dt$. There are $\binom{dt}{2j}(2j-1)!!$ collections of $j$ disjoint pairs among the $dt$ half-edges of $U$. Each prescribed collection appears in the random pairing with probability $\prod_{i=0}^{j-1}(dn-2i-1)^{-1}$. A union bound therefore gives
    \begin{align*}
        \Pb\Big( |E(G[U])|\ge j\Big) &\le \binom{dt}{2j}(2j-1)!! \prod_{i=0}^{j-1}\frac{1}{dn-2i-1}\\
        &\le \frac{(dt)^{2j}}{2^j j!} \left(\frac{2}{dn}\right)^j \le \left(\frac{edt^2}{jn}\right)^j \le \left(\frac{edt}{n} \right)^j.
    \end{align*}
    Here we used $2j\le dt\le\delta dn<dn/2$, $j!\ge(j/e)^j$, and $j\ge t$. Since $\e d\,\delta<1$ and $j>(1+\epsilon)t$, a union bound over the vertex sets gives
    \begin{align*}
        \Pb\Big( \exists U\subseteq V: |U|=t, |E(G[U])|>(1+\epsilon)t \Big) \le\binom nt\left(\frac{edt}{n}\right)^j \le\left[ e(ed)^{1+\epsilon} \left(\frac tn\right)^\epsilon\right]^t.
    \end{align*}
    Write $C=e(ed)^{1+\epsilon}$. Our choice of $\delta$ gives
    \[
        C\delta^\epsilon =\frac{e(ed)^{1+\epsilon}}{(100d)^2} \le \frac{e^3}{10^4}<\frac12.
    \]
    For $t\le\sqrt n$, we also have $C(t/n)^\epsilon\le Cn^{-\epsilon/2}$. Summing over all sizes, we obtain
    \[
        \Pb\Big( \exists U: 1\le|U|\le\delta n, |E(G[U])|>(1+\epsilon)|U| \Big) \le \sum_{t=1}^{\lfloor\sqrt n\rfloor}(Cn^{-\epsilon/2})^t +\sum_{t>\sqrt n} 2^{-t} =o(1). \qedhere
    \]
\end{proof}
Finally, we prove that Item~(5) holds for $G \sim\mathcal G_d(n)$ with probability $1-o(1)$, as incorporated in the following lemma.
\begin{lemma}{\label{lem:small-ball-removal-degree}}
    The following holds with probability $1-o(1)$ if $G\sim\mathcal G_d(n)$ with $d \ge 3$. Recall \eqref{eq:def-radii-r-R}, \eqref{eq:def-local-ball}, and \eqref{eq:structural-parameters}. For all $v \in V$ and for all $H \subseteq F_{v,r}(E)$ such that $H$ is connected and $|V(H)| \ge \log n$, we have
        \[
            \big| V(H) \cap S_{v,r} \big| \le \frac{\epsilon}{d}|V(H)|, \qquad \sum_{x \in V(H)} \big( d-\deg_{F_{v,r}(E)}(x) \big) \le \epsilon |V(H)|.
        \]
\end{lemma}
\begin{proof}
    We work on the event of Item~(1). Fix $v\in V$, and let $K$ be the graph obtained from $G[B_{v,R}]$ by removing all edges of $G[B_{v,r}]$. Let $\mathcal C$ be the collection of components of $K$ intersecting $S_{v,r}$. Since $G[B_{v,r}]$ is connected, we may take a spanning tree of $G[B_{v,r}]$ and one of each $C\in\mathcal C$. These trees have disjoint edge sets, and the tree of $C$ meets the inner ball exactly at $V(C)\cap S_{v,r}$. Adding that tree increases the excess by $|V(C)\cap S_{v,r}|-1$. Thus
    \[
        \sum_{C\in\mathcal C}\bigl(|V(C)\cap S_{v,r}|-1\bigr)
        \le  \tau(G[B_{v,R}])\le1.
    \]
    Thus each component of $K$ contains at most one vertex in $S_{v,r}$, except that at most one component may contain two. Now let $H\subseteq F_{v,r}(E)$ be connected, and write $t=|V(H)|\ge\log n$ and $k=|V(H)\cap S_{v,r}|$. By discarding at most one vertex from $V(H)\cap S_{v,r}$, we obtain a set $W$ whose vertices lie in distinct components of $K$. Put $L=\lceil 2d/\epsilon \rceil$. For all sufficiently large $n$, we have $r+L<R$ and $\log n\ge2(L+1)$. Thus any two vertices of $W$ have distance greater than $2L$ in $F_{v,r}(E)$. The balls of radius $L$ in $H$ centered at vertices of $W$ are therefore pairwise disjoint. Each contains at least $L+1$ vertices, since $H$ is connected and $t\ge L+1$. Consequently,
    \[
        |W|(L+1)\le t, \qquad k\le |W|+1\le\frac{t}{L+1}+1 \le \frac{\epsilon}{d} t.
    \]
    This proves the first inequality, simultaneously for all $v$ and $H$.

    Finally, every vertex of $B_{v,r}\setminus S_{v,r}$ is isolated in $F_{v,r}(E)$. Since $H$ is connected and $t\ge2$, it contains no such vertex. Vertices outside $B_{v,r}$ retain degree $d$, therefore
    \[
        \sum_{x\in V(H)}\Big( d-\mathsf{deg}_{F_{v,r}(E)}(x) \Big) \le d|V(H)\cap S_{v,r}| \le \epsilon|V(H)|. \qedhere
    \]
\end{proof}

\section{Random-cluster model on trees}{\label{sec:random-cluster-model-on-tree}}

We record in this section some preliminary inputs for the random-cluster model on trees and tree-like graphs. We first define the random-cluster model with boundary conditions. Let $\Lambda \subseteq V(G)$ be an arbitrary subset of vertices. A boundary condition on $\Lambda$ is a partition $\xi=(P_1,\ldots,P_k)$ of $\Lambda$. Vertices belonging to the same part of the partition are said to be wired together. As such, the random-cluster distribution on $G$ with parameters $p \in (0,1)$ and $q>0$ and boundary condition $\xi$ on $\Lambda\subseteq V(G)$ is defined as
\begin{equation}{\label{eq:def-random-cluster-model-boundary-condition}}
    \mu_{G,p,q}^{\xi}(A) \propto p^{|A|} (1-p)^{|E \setminus A|} q^{c(A,\xi)},
\end{equation}
where $c(A,\xi)=c_G(A,\xi)$ denotes the number of connected components in the subgraph $(V,A)$ when taking into account the wirings in $\xi$. We write $\xi\preceq\xi'$ if every block of $\xi$ is contained in a block of $\xi\prime$. For $q\ge1$, coarsening the boundary partition increases the random-cluster measure in stochastic order. The spatial Markov property identifies the conditional interior law with the measure whose boundary partition is induced by the exterior configuration \cite{Grimmett06}. 

\subsection{The cavity recursion}{\label{subsec:cavity-recursion}}

For $D \in \{ d-1,d \}$, let $\mathbb T_{h,D}$ be the rooted tree of height $h$ with $D$ children at the root and $d-1$ children at every other internal vertex. Unless otherwise specified, boundary conditions on a finite tree $\mathbb T$ are partitions of $\Lambda=\mathsf L(\mathbb T)$, the set of leaves of $\mathbb T$. For any tree $\mathbb T$ and $U \subseteq \mathsf L(\mathbb T)$, let $\xi_U$ identify the vertices of $U$ and leave all other leaves as singletons. For any $u \in \mathbb T$, denote by $\mathbb T_u$ the descendant subtree of $u$ in $\mathbb T$. We will use the cavity tree recursion and its fixed-point analysis from \cite[Section~2]{BGPZ26}. Put
\begin{equation}{\label{eq:cavity-function}}
    \Phi(x) = \frac{ x+(q-1)(1-p) }{ (1-p)x+p+(q-1)(1-p) }, \qquad F(x)=\Phi(x)^{d-1}.
\end{equation}
Represent the wired block by a ghost vertex and omit its constant component
factor. Fixing a single-component boundary condition $\xi=\xi_U$, let $\mathcal Z_{\mathbb T}^{\xi,0}$ denote the contribution to the partition function of the random-cluster configurations on $\mathbb T$ under the boundary condition $\xi$ in which $\rho(\mathbb T)$ is not connected to the wired component $U$; analogously define $\mathcal Z_{\mathbb T}^{\xi,1}$ for the total weight of the configurations that connect $\rho(\mathbb T)$ to $U$. For any vertex $u \in V(\mathbb T)$, let $\xi(u)$ be the single-component boundary condition induced on $\mathsf L(\mathbb T_u)$ by $\xi$. For any finite rooted $\mathbb T$ and any single-component boundary condition $\xi=\xi_U$ where $U \subseteq \mathsf L(\mathbb T)$, define the message function $\mathbf X:V(\mathbb T) \to [1,\infty]$ by
\begin{equation}{\label{eq:message-function}}
    \mathbf X_{\mathbb T}^{\xi}(u) = q \frac{ \mathcal Z_{\mathbb T_u}^{\xi(u),1} }{ \mathcal Z_{\mathbb T_u}^{\xi(u),0} } +1.
\end{equation}
We will record the following properties proved in \cite[Section~2]{BGPZ26}.
\begin{lemma}{\label{lem:cavity-fixed-point}}
    When $q>2$ and $p>\ps$, the following conditions hold.
    \begin{itemize}
        \item[(1)] For any finite rooted tree $\mathbb T$ and any single-component boundary condition $\xi=\xi_U$ for some $U \subseteq \mathsf L(\mathbb T)$, the leaf messages are $\infty$ on $U$ and $1$ on $\mathsf L(\mathbb T)\setminus U$. At every nonleaf vertex $u$, the message function satisfies
        \[
            \mathbf X_{\mathbb T}^{\xi}(u) = \prod_{ w:u \to w } \Phi(\mathbf X_{\mathbb T}^{\xi}(w)).
        \]
        Here $\Phi(\infty)=(1-p)^{-1}$. In particular, when $\mathbb T\in\{\mathbb T_{h,d-1},\mathbb T_{h,d}\}$, every inward message from $u \in V(\mathbb T) \setminus (\mathsf L(\mathbb T) \cup \{ \rho(\mathbb T) \})$ satisfies $1 \le \mathbf X_{\mathbb T}^{\xi}(u) \le (1-p)^{-(d-1)}$.
        \item[(2)] There exists a unique $y_{\star}>1$ with $F(y_\star)=y_\star$. In addition,
        \[
            0<\beta_\star:=F'(y_\star)<1, \qquad F^{\circ k}(x) \to y_\star \mbox{ for all } x>1.
        \]
        The convergence is monotone from either side of $y_\star$.
    \end{itemize}
\end{lemma}
We will also use the root-edge marginal in terms of its two incoming messages. For each root edge $e=(u,v) \in E_\rho(\mathbb T)$, cutting the edge $e$ in $\mathbb T$ gives two subtrees rooted at $u$ and $v$, respectively, with the same ghost vertex representing the wired block. If their messages are $x$ and $y$, summing over the state of the cut edge gives
\begin{equation}{\label{eq:root-edge-marginal}}
\begin{aligned}
    R(x,y) &= \mu^{\xi}_{\mathbb T}( (u,v) \mbox{ open} ) \\
    &= \frac{ p(xy+q-1) }{ xy+(q-1)(1-p)(x+y) + (q-1)(p+(q-1)(1-p)) }.
\end{aligned}
\end{equation}
We next record a sharper estimate when comparing two boundary conditions $\xi_U,\xi_{U'}$ with $U \subseteq U'$. Motivated by \cite{LS26}, we first introduce the following definition of a dense boundary condition.
\begin{defn}{\label{def:dense-boundary-condition}}
    Put $\ell=\ell_h=2\lceil \log_{d-1}h \rceil$. For $h$ large enough that $\ell_h<h$ and $\mathbb T \in \{ \mathbb T_{h,d-1}, \mathbb T_{h,d} \}$, a set $W \subseteq \mathsf L(\mathbb T)$ is $\kappa$-dense if every descendant subtree of height $\ell$ contains at least $\kappa (d-1)^{\ell}$ leaves in $W$.
\end{defn}
\begin{lemma}\label{lem:dense-boundary-comparison}
    Suppose $q>2$ and $p>\ps$, and put $\sigma=\frac{1+\beta_\star}{2}$. Consider $0<\kappa\le1$, $\mathbb T\in \{ \mathbb T_{h,d-1}, \mathbb T_{h,d} \}$, and
    \[
        W \subseteq U \subseteq U' \subseteq \mathsf L(\mathbb T), \qquad W \mbox{ is } \kappa\text{-dense}.
    \]
    Then for every $e \in E_\rho(\mathbb T)$,
    \[
        0 \le \mu_{\mathbb T}^{\xi_{U'}}(e \mbox{ open}) - \mu_{\mathbb T}^{\xi_{U}}(e \mbox{ open}) \le O_{\kappa}(1) \cdot |U' \setminus U| \left( \frac{\sigma}{d-1} \right)^h.
    \]
\end{lemma}
\begin{proof}
    It suffices to consider $h$ sufficiently large depending on $\kappa$, since the remaining finitely many heights can be absorbed into the implicit constant. Put $M=(1-p)^{-(d-1)}$. The logarithmic derivative of $\Phi$ is
    \[
        \psi(x):=\frac{x\Phi'(x)}{\Phi(x)} =\frac{p(p+q(1-p))x}{(x+(q-1)(1-p))((1-p)x+p+(q-1)(1-p))}.
    \]
    The expression above gives $0\le\psi(x)\le1$ for $x\ge1$. Moreover, since $F(y_\star)=y_\star$, Lemma~\ref{lem:cavity-fixed-point} gives
    \[
        (d-1)\psi(y_\star)=F'(y_\star)=\beta_\star<1.
    \]
    We first show that density keeps the messages at height $\ell$ above the leaves uniformly above $1$. Since $(d-1)\widehat p>1$, there is a constant $\varphi\in(0,1)$ such that
    \[
        1-(1-\widehat p\varphi)^{(d-1)} \ge \varphi.
    \]
    Under independent Bernoulli($\widehat p$) percolation on a complete
    $(d-1)$-ary tree, the probability $s_j$ that its root reaches generation $j$ satisfies $s_0=1$ and $s_{j+1}=1-(1-\widehat p s_j)^{(d-1)}$. Hence $s_j\ge\varphi$ for every $j$. Conditional on reaching generation $\ell$, choose uniformly one of the reached leaves. By symmetry, this chosen leaf is uniform over all $(d-1)^\ell$ leaves. Thus the probability of reaching any prescribed set containing at least $\kappa (d-1)^\ell$ leaves is at least $\kappa\varphi$. The random-cluster measure stochastically dominates independent Bernoulli($\widehat p$) percolation. Apply this domination in each descendant subtree of height $\ell$, with the boundary condition induced by any wired block containing $W$. The root of that subtree connects to the ghost with probability at least $\kappa\varphi$. If its message is $x$, then \eqref{eq:message-function} expresses this probability as
    \[
        \frac{\mathcal Z^1}{\mathcal Z^0+\mathcal Z^1} =\frac{x-1}{x+q-1}.
    \]
    Consequently every such message belongs to $[x_-,M]$, where
    \[
        x_-:=1+\min\left\{ \frac{q\kappa\varphi}{1-\kappa\varphi}, \frac{y_\star-1}{2} \right\}>1.
    \]
    Notice that $x_-<y_\star<M$. Fix $a=\frac{\beta_\star+\sigma}{2(d-1)}$. By monotonicity of $\Phi$, a message at height $\ell+j$ above the leaves, at a vertex other than the root, lies in
    \[
        I_j:=[F^{\circ j}(x_-),F^{\circ j}(M)].
    \]
    Lemma~\ref{lem:cavity-fixed-point} shows that these intervals are nested and shrink to $\{y_\star\}$. Since $\psi(y_\star)<a$, we can choose an integer $K=K(\kappa)$ such that $\psi(x)\le a$ for every $x\in I_K$. It follows that, under every boundary condition whose wired block contains $W$, all messages at depths $1,\ldots,h-\ell-K$ belong to $I_K$. 

    Now add a single leaf $z\notin U$ to the wired block. Its message changes from $1$ to $\infty$, so the logarithm of its parent's message increases by exactly $\log(1/(1-p))$. Only messages on the path from $z$ to the root can change. At each subsequent ancestor, the recursion in Lemma~\ref{lem:cavity-fixed-point} and the mean value theorem in logarithmic coordinates give
    \[
        \log\Phi(x')-\log\Phi(x)
        =\psi(\widetilde x)(\log x'-\log x)
        \qquad\text{for some }\widetilde x\in[x,x'],
    \]
    where $x\le x'$ are the two messages from the affected child. Each factor is at most $1$. For affected children at depths $2,\ldots,h-\ell-K$, both messages lie in $I_K$, so the corresponding $h-\ell-K-1$ factors are at most $a$. The change in the logarithm of the message reaching the child of the root is therefore at most $\log(1/(1-p)) \cdot a^{h-\ell-K-1}$.

    Fix a root edge $e$. Of its two cut-branch messages in \eqref{eq:root-edge-marginal}, only the one on the side containing $z$ changes. If $z$ lies in a different root branch from $e$, one additional factor, bounded by $1$, transfers the change to the root-side message. Both cut-branch messages lie in $[1,M]$: the child side has $(d-1)$ children, while the root side has $\le (d-1)$ children. The formula for $R(x,y)$ is smooth on this compact square, and hence its partial derivatives with respect to $\log x$ and $\log y$ are bounded by a constant when $d,p,q$ are fixed. We thus obtain
    \[
        0\le
        \mu_{\mathbb T}^{\xi_{U\cup\{z\}}}(e\mbox{ open})
        -\mu_{\mathbb T}^{\xi_U}(e\mbox{ open})
        \le O(1) \cdot a^{h-\ell-K-1} \le O_\kappa(1) \cdot (\sigma/(d-1))^h.
    \]
    Add the leaves of $U'\setminus U$ one at a time and sum these bounds. Every intermediate wired block still contains $W$, so the same constants apply at each step. This proves the lemma.
\end{proof}

Lemma~\ref{lem:dense-boundary-comparison} treats single-block boundary conditions. We next allow additional common boundary blocks. The estimate remains proportional to the number of boundary vertices on which the two partitions may differ.
\begin{lemma}\label{lem:multi-block-boundary-comparison}
    Suppose $q>2$ and $p>\ps$, and choose $\sigma\in(0,1)$ as in Lemma~\ref{lem:dense-boundary-comparison}. Consider $0<\kappa\le1$ and $\mathbb T\in \{ \mathbb T_{h,d-1}, \mathbb T_{h,d} \}$.
    %For every $0<\kappa\le1$, there are constants $C_0<\infty$ and $h_0>0$, depending only on $\kappa$, with the following property. Let $h\ge h_0$, $\mathbb T=\mathbb T_{h,D}$ with $D\in\{d-1,d\}$. 
    Suppose $S \subseteq \mathsf L(\mathbb T)$ and let $\xi,\xi'$ have the same restricted partition $\pi$ on $\mathsf L(\mathbb T) \setminus S$. Suppose a block of $\pi$ contains a $\kappa$-dense set $W$, and that $\pi$ has $J$ other non-singleton blocks. Then
    \[
        \TV\left( \mu^{\xi}_{\mathbb T}\big|_{E_\rho(\mathbb T)},  \mu^{\xi'}_{\mathbb T}\big|_{E_\rho(\mathbb T)} \right) \le O_\kappa(1)\cdot |S| \big(1+J(h+1)^{O(1)}\big) \left( \frac{\sigma}{d-1} \right)^h.
    \]
\end{lemma}
\begin{proof}
    If $S=\emptyset$, then $\xi=\xi'$ and the assertion is immediate.
    Henceforth put $Q=|S|\ge1$ and $L=\mathsf L(\mathbb T)$. Let $U$ be the
    block of $\pi$ containing $W$, and let $B_1,\ldots,B_J$ be its other
    non-singleton blocks. Write $\xi_1=\xi$ and $\xi_2=\xi'$. Put
    \[
        N_h=\left\lceil \left( 1+\log\left( \frac{d-1}{\sigma} \right) \right) (d-1)^\ell \widehat p^{-2\ell} h \right\rceil = (h+1)^{O(1)}.
    \]
    We call $B_j$ large if $|B_j|\ge N_h$, and small otherwise. Let $\zeta$ be a boundary partition with a block containing $W$, and let $B$ be a set of leaves in another block. At least $|B|/(d-1)^\ell$ of the bottom descendant subtrees of height $\ell$ meet $B$. Each such subtree contains a leaf of $W$ as $W$ is dense. In each such subtree, choose a leaf of $B$, a leaf of $W$, and the tree path joining them. These paths are edge-disjoint and have length at most $2\ell$. Since every random-cluster measure with the given parameters stochastically dominates a bond percolation with parameter $\widehat p$, we obtain
    \begin{equation}\label{eq:large-boundary-block-connection}
        \mu_{\mathbb T}^{\zeta}(B\not\leftrightarrow W) \le\exp\left(-\frac{|B|\widehat p^{2\ell}}{(d-1)^\ell}\right) \le \left( \frac{\sigma}{d-1} \right)^h.
    \end{equation}
    Suppose $\zeta^+$ is obtained by merging some blocks $\{ B_1,\ldots,B_m \}$ of $\zeta$ into its block containing $W$. The definition of the random-cluster measure gives
    \[
        \mu_{\mathbb T}^{\zeta^+}(A) =\mu_{\mathbb T}^{\zeta}(A) \frac{R(A)}{\mathbb E_{\mu_{\mathbb T}^{\zeta}}R}, \qquad \Delta(A)=c(A,\zeta)-c(A,\zeta^+)\ge0, \qquad R(A)=q^{-\Delta(A)}\le1.
    \]
    Since $\mathbb E_{\mu_{\mathbb T}^{\zeta}}R\le1$, both probabilities at
    $A$ are at least $\mu_{\mathbb T}^{\zeta}(A)R(A)$. Moreover, $R(A)=1$
    whenever all the merged blocks are already connected to $W$. Therefore
    \begin{equation}\label{eq:boundary-block-merging-cost}
        \TV(\mu_{\mathbb T}^{\zeta^+},\mu_{\mathbb T}^{\zeta})
        \le1-\mathbb E_{\mu_{\mathbb T}^{\zeta}}R \le\mu_{\mathbb T}^{\zeta}\Big( \exists B \in \{ B_1,\ldots,B_m \}, B\not\leftrightarrow W \Big).
    \end{equation}

    For $i=1,2$, form $\xi_i^+$ by merging the full blocks of $\xi_i$
    containing the large $B_j$ into its full block containing $U$. Distinct
    blocks of $\pi$ lie in distinct full blocks of $\xi_i$, since otherwise
    the restriction of $\xi_i$ would not equal $\pi$. Consequently, the
    restrictions of $\xi_1^+$ and $\xi_2^+$ to $\mathsf L(\mathbb T) \setminus S$ still agree.
    Their common restricted partition has one large block
    \[
        U^+=U\cup\left(\cup_{j:\,|B_j|\ge N_h}B_j \right),
    \]
    and its other non-singleton blocks are precisely the small $B_j$.
    By \eqref{eq:large-boundary-block-connection}, \eqref{eq:boundary-block-merging-cost}, and a union bound,
    \begin{equation}{\label{eq:TV-xi-xi-plus}}
        \TV(\mu_{\mathbb T}^{\xi_i^+},\mu_{\mathbb T}^{\xi_i}) \le J\left(\frac{\sigma}{d-1}\right)^h, \qquad i=1,2.
    \end{equation}

    We now place the two modified partitions between single-block boundary
    conditions. Put
    \[
        F_i=\Big\{ v \in\mathsf L(\mathbb T): v \notin S \cup U^+ \cup \cup_{|B_j|<N_h} B_j, v \mbox{ belongs to a non-singleton block of } \xi_i^+ \Big\}, \quad i=1,2.
    \]
    For each leaf $v \in F_i$, its full block contains at least one vertex of $S$ and no other vertex outside $S$. Distinct such blocks use disjoint subsets of $S$, so $|F_i|\le |S|$. Put
    \[
        B_*=S\cup F_1\cup F_2 \cup \left( \cup_{j:|B_j|<N_h} B_j \right).
    \]
    Then $B_*\cap U^+=\emptyset$ and $|B_*|\le 3|S|+JN_h$. Clearly 
    \[
        \xi_{U^+}\preceq\xi_i^+\preceq\xi_{U^+\cup B_*}, \qquad i=1,2.
    \]
    The coupling inequality and a union bound over the at most $d$ root edges therefore give
    \begin{align*}
        \TV\left( \mu_{\mathbb T}^{\xi_1^+}\big|_{E_\rho(\mathbb T)}, \mu_{\mathbb T}^{\xi_2^+}\big|_{E_\rho(\mathbb T)} \right) &\le\sum_{e\in E_\rho(\mathbb T)}
        \left[ \mu_{\mathbb T}^{\xi_{U^+\cup B_*}}(e\mbox{ open}) -\mu_{\mathbb T}^{\xi_{U^+}}(e\mbox{ open}) \right]\\
        &\le O_\kappa(1) (3|S|+JN_h) \left(\frac{\sigma}{d-1}\right)^h.
    \end{align*}
    The last inequality follows from Lemma~\ref{lem:dense-boundary-comparison}, applied with the same dense set $W\subseteq U^+$.

    Combined with \eqref{eq:TV-xi-xi-plus} and the triangle inequality,
    \begin{align*}
        \TV\left( \mu_{\mathbb T}^{\xi}\big|_{E_\rho(\mathbb T)}, \mu_{\mathbb T}^{\xi'}\big|_{E_\rho(\mathbb T)} \right) &\le \Big[ 2J+O_\kappa(1)(3|S|+JN_h) \Big] \left(\frac{\sigma}{d-1}\right)^h \\
        &\le O_\kappa(1)\cdot |S| \Big( 1+J(h+1)^{O_\kappa(1)} \Big) \left(\frac{\sigma}{d-1}\right)^h. \qedhere
    \end{align*}
\end{proof}

\subsection{Mixing on graphs with small excess}{\label{subsec:mixing-small-excess-graphs}}

To handle neighborhoods containing a cycle, we use the weak spatial mixing
estimate of \cite[Theorem~5.1]{BGPZ26}. We first specify its boundary and
graph hypotheses in our notation. The following definition is the coupling
formulation of \cite[Definition~4.4]{BGPZ26}.
\begin{defn}
    A distribution $\mathcal L$ over boundary conditions $\xi$ on $\Lambda$ is $(\theta,Q)$-wired if one can couple $\xi\sim\mathcal L$ to a random subset $Z\subseteq\Lambda$, obtained by including each vertex independently with probability $\theta$, and an exceptional set $U\subseteq\Lambda$ such that $|U|\le Q$ almost surely. We require a set $W\subseteq\Lambda\setminus U$ with $Z\setminus U\subseteq W$ such that the restriction of $\xi$ to $\Lambda\setminus U$ consists of the block $W$ and singleton blocks for all remaining vertices. The vertices of $U$ may be wired arbitrarily; both $U$ and $W$ may depend on $Z$ and $\xi$.
\end{defn}
\begin{defn}[Almost regular unicyclic graphs] \label{def:almost-regular-unicyclic}
    Fix $d\ge3$ and an integer $h\ge1$. Let $\mathcal U_d(h)$ be the class
    of finite connected simple graphs $H$ with a root $\rho$ and a
    BFS tree $T$ rooted at $\rho$, satisfying the
    following conditions:
    \begin{itemize}
        \item[(1)] Every leaf of $T$ has depth $h$.
        \item[(2)] $|E(H)\setminus E(T)|\le1$, or equivalently,
        $\tau(H)\le1$.
        \item[(3)] $\deg_H(\rho)\in\{d-1,d\}$, and
        $\deg_H(v)=d$ for every
        $v\in V(H)\setminus(\{\rho\}\cup\mathsf L(T))$.
    \end{itemize}
    The boundary of such a graph is
    \[
        \Lambda_H:=\mathsf L(T)
        =\{v\in V(H):\dist_H(\rho,v)=h\}.
    \]
    We allow the tree case $\tau(H)=0$; in this case
    $H\in\{\mathbb T_{h,d-1},\mathbb T_{h,d}\}$, up to rooted isomorphism.
    Write $E_\rho(H)$ for the edges of $H$ incident to $\rho$, and write
    $\mathrm w$ for the boundary condition identifying all vertices of
    $\Lambda_H$.
\end{defn}
For an edge $e=\{u,v\}\in E(H)$, put
\[
    s_e:=\dist_H(e,\Lambda_H)
    =\min\{\dist_H(u,\Lambda_H),\dist_H(v,\Lambda_H)\}.
\]
In particular, $s_e=h-1$ for every $e\in E_\rho(H)$. We record the
edge-event conclusion of the cited spatial mixing theorem, replacing its
asymptotic rate by a fixed rate strictly larger than $\beta_\star$.
\begin{lemma}\label{lem:unicyclic-spatial-mixing}
    Fix $d\ge3$, $q>2$, $\ps<p<1$, $\theta\in(0,1]$, and an integer $Q\ge0$. Let $\beta_\star=F'(y_\star)$ be as in Lemma~\ref{lem:cavity-fixed-point}, and fix $\sigma\in(\beta_\star,1)$. There are constants $c,C>0$ and an integer $h_0\ge1$, depending only on $d,p,q,\theta,Q,\sigma$, such that the following holds. Let $H\in\mathcal U_d(h)$, and let $\mathcal L$ be a $(\theta,Q)$-wired distribution on partitions of $\Lambda_H$. For every $e\in E(H)$ with $s_e\ge h_0$,
    \begin{equation}\label{eq:unicyclic-edge-response}
        \Pb_{\xi\sim\mathcal L}\left( 0\le\mu_H^{\mathrm w}(e\mbox{ open}) -\mu_H^\xi(e\mbox{ open}) \le C\sigma^{s_e/2} \right)  \ge1-\exp\{-c(d-1)^{\sqrt{s_e}}\}.
    \end{equation}
    Consequently, if $h-1\ge h_0$, with probability at least $1-d\exp\{-c(d-1)^{\sqrt{h-1}}\}$ over $\xi\sim\mathcal L$,
    \[
        \TV\left( \mu_H^\xi\big|_{E_\rho(H)}, \mu_H^{\mathrm w}\big|_{E_\rho(H)} \right) \le dC\sigma^{(h-1)/2}.
    \]
\end{lemma}
\begin{proof}
    The first assertion is the edge-event conclusion of \cite[Theorem~5.1]{BGPZ26}. For the second assertion, apply \eqref{eq:unicyclic-edge-response} to each of the at most $d$ root edges and take a union bound. On the resulting event, a monotone coupling of $\mu_H^\xi \preceq \mu_H^{\mathrm w}$ and the coupling inequality give
    \[
        \TV\left( \mu_H^\xi\big|_{E_\rho(H)}, \mu_H^{\mathrm w}\big|_{E_\rho(H)} \right)
        \le\sum_{e\in E_\rho(H)} \bigl[\mu_H^{\mathrm w}(e\mbox{ open}) -\mu_H^\xi(e\mbox{ open})\bigr] \le dC\sigma^{(h-1)/2}. \qedhere
    \]
\end{proof}

Finally, we prove a spectral gap estimate for the random-cluster dynamics on graphs with small excess. 
\begin{lemma}{\label{lem:spec-gap-small-excess}}
    Fix $p\in(0,1)$ and $q\ge1$. Let $H$ be a finite connected graph with $M:=|E(H)|\ge1$ edges, allowing loops and parallel edges, and let $\xi$ be a boundary condition with $J$ non-singleton blocks. Suppose $H$ has a rooted spanning tree of height at most $h$, where $h\ge0$ is an integer. Put
    \[
        R = \frac{ |E(H)|^2 q^{ 2h+4\tau(H)+2J } }{ \min\{ \widehat{p},1-p \} }.
    \]
    The discrete-time, uniform-edge heat-bath kernel $P^{\xi}_H$
    satisfies, for every $0<\delta<1$,
    \[
        \mathsf{gap}(P^{\xi}_H) \ge R^{-1}, \qquad T_\mix(P^{\xi}_H,\delta) \le R\big( |E(H)|\log(q/\min\{ p,1-p \}) + \log(1/\delta) \big).
    \]
    Here $T_\mix(P^{\xi}_H,\delta)$ denotes the mixing time to
    total-variation distance at most $\delta$ from a worst-case initial
    configuration.
\end{lemma}
\begin{proof}
    We use the canonical-path bound of \cite[Corollary~13.21]{LPW17}
    with an edge ordering obtained from the spanning tree. Write
    $\mu=\mu_H^\xi$ and $P=P_H^\xi$, and let $\overline H=H/\xi$
    be the graph obtained by identifying the vertices in each block
    of $\xi$. Retain all edge labels, including those of loops and
    parallel edges, so that configurations on $H$ and $\overline H$
    are indexed by the same subsets of $E(H)$. For each configuration
    $A$, the number $c(A,\xi)$ is then the component count of
    $(V(\overline H),A)$.

    Choose a rooted spanning tree $T$ of height at most $h$. List
    its edges in depth-first order, placing each edge to a child
    immediately before the edges in that child's descendant subtree.
    Append the $\tau(H)$ edges outside $T$ in any order. For the
    resulting ordering $e_1,\ldots,e_M$, set
    \[
        F_i=\{e_1,\ldots,e_i\},\qquad
        K_i=E(H)\setminus F_i\qquad(0\le i\le M),
    \]
    and let $C_i$ be the set of vertices of $\overline H$ incident
    to edges in both $F_i$ and $K_i$. We claim that
    \begin{equation}\label{eq:small-excess-active-set}
        |C_i|\le w:=h+2\tau(H)+J\qquad(0\le i\le M).
    \end{equation}
    During the tree part of the ordering, the vertices of $H$
    incident to tree edges on both sides of the cut lie on the
    current ancestor path. Each has depth less than $h$, so there
    are at most $h$ such vertices. Every other vertex incident to
    both sides is an endpoint of a non-tree edge, contributing at
    most $2\tau(H)$ vertices. During the appended part, every vertex
    incident to both sides is an endpoint of a remaining non-tree
    edge. An unmerged vertex belongs to $C_i$ only if it meets
    both sides in $H$. The $J$ nonsingleton boundary blocks contribute at
    most $J$ further vertices, proving
    \eqref{eq:small-excess-active-set}.

    This bound controls the random-cluster weights when two
    configurations exchange their edges across a cut. Fix $i$ and
    configurations $A,B$, and define
    \[
        U=(B\cap F_i)\cup(A\cap K_i),\qquad
        W=(A\cap F_i)\cup(B\cap K_i).
    \]
    We will show that
    \begin{equation}\label{eq:small-excess-weight-comparison}
        \mu(A)\mu(B)\le q^{2w}\mu(U)\mu(W).
    \end{equation}
    Split each vertex of $C_i$ into two copies, assigning its
    incident edges in $F_i$ to one copy and those in $K_i$ to the
    other. Keep every other vertex on the side containing its
    incident edges, and assign any isolated vertices of $\overline H$
    to the $F_i$ side. The vertex sets of the two resulting graphs
    depend only on the cut. Write $k_i^-(A\cap F_i)$ and
    $k_i^+(A\cap K_i)$ for their component counts, including vertices
    isolated by the choice of $A$. Identifying each pair of copies
    recovers $(V(\overline H),A)$ and decreases the component count
    by at most one per identification. Thus
    \[
        c(A,\xi)=k_i^-(A\cap F_i)+k_i^+(A\cap K_i)-\delta_i(A),
        \qquad 0\le\delta_i(A)\le |C_i|.
    \]
    The component counts on each side cancel when comparing
    $(A,B)$ with $(U,W)$, giving
    \begin{align*}
        c(A,\xi)+c(B,\xi)-c(U,\xi)-c(W,\xi)
        &=\delta_i(U)+\delta_i(W)-\delta_i(A)-\delta_i(B)\\
        &\le2w.
    \end{align*}
    Since $|A|+|B|=|U|+|W|$, the product edge weights and normalizing
    constants also cancel. The assumption $q\ge1$ now gives
    \eqref{eq:small-excess-weight-comparison}.

    For each ordered pair $(A,B)$, let $\gamma_{A,B}$ run from $A$
    to $B$ by changing the edges of $A\mathbin{\triangle}B$ one at
    a time in the prescribed order. These paths have length at most
    $M$. Fix a directed transition $X\to X'$ that changes $e_r$.
    If $\gamma_{A,B}$ uses this transition, then
    \[
        X=(B\cap F_{r-1})\cup(A\cap K_{r-1}).
    \]
    Associate to its endpoints the configuration
    $Z=(A\cap F_{r-1})\cup(B\cap K_{r-1})$. For the fixed
    transition, $Z$ determines both endpoints through
    \[
        A=(Z\cap F_{r-1})\cup(X\cap K_{r-1}),\qquad
        B=(X\cap F_{r-1})\cup(Z\cap K_{r-1}).
    \]
    This injectivity and \eqref{eq:small-excess-weight-comparison}
    at the cut $r-1$ imply
    \[
        \sum_{(A,B):\,(X,X')\in\gamma_{A,B}}
            \mu(A)\mu(B)|\gamma_{A,B}|
        \le Mq^{2w}\mu(X)\sum_{Z\subseteq E(H)}\mu(Z)
        =Mq^{2w}\mu(X).
    \]
    Conditional on selecting an edge, its opening probability is
    either $\widehat p$ or $p$, and equals $p$ for loops.
    Every opening therefore has probability at least $\widehat p$,
    and every closing has probability at least $1-p$. Since the
    edge is chosen uniformly, every allowed flip satisfies
    $P(X,X')\ge\min\{\widehat p,1-p\}/M$. Hence the congestion is
    at most
    \[
        \max_{\substack{X\ne X'\\P(X,X')>0}}
        \frac{1}{\mu(X)P(X,X')}
        \sum_{(A,B):\,(X,X')\in\gamma_{A,B}}
            \mu(A)\mu(B)|\gamma_{A,B}|
        \le\frac{M^2q^{2w}}{\min\{\widehat p,1-p\}}=R.
    \]
    The canonical-path bound gives $\mathsf{gap}(P)\ge R^{-1}$.

    To pass to the mixing time, put $u_0=\min\{p,1-p\}$. We first
    bound the minimum stationary mass by
    \begin{equation}\label{eq:small-excess-minimum-mass}
        \min_{A\subseteq E(H)}\mu(A)\ge(u_0/q)^M.
    \end{equation}
    Indeed, each unnormalized configuration weight is at least
    $u_0^Mq$, and the partition function is at most
    $q^{|V(\overline H)|}$ because the product edge weights sum to
    one. Connectivity gives $|V(\overline H)|-1\le M$, which yields
    \eqref{eq:small-excess-minimum-mass}.

    Each single-edge heat-bath operator is an orthogonal
    conditional-expectation projection in $L^2(\mu)$, so their
    average $P$ has nonnegative spectrum. The spectral estimate in
    \cite[proof of Theorem~12.4]{LPW17}, together with
    \eqref{eq:small-excess-minimum-mass}, therefore gives, for every
    integer $t\ge0$,
    \[
        \max_A\TV\bigl(P^t(A,\cdot),\mu\bigr)
        \le\tfrac12\e^{-t/R}(q/u_0)^{M/2}.
    \]
    Set
    \[
        t_\delta=\left\lfloor
            R\bigl[M\log(q/u_0)+\log(1/\delta)\bigr]
        \right\rfloor.
    \]
    Since $\min\{\widehat p,1-p\}\le1/2$, we have $R\ge2$.
    Also $M\ge1$ and $q/u_0\ge2$, so
    \[
        \max_A\TV\bigl(P^{t_\delta}(A,\cdot),\mu\bigr)
        \le\frac{\delta}{2}\e^{1/R}(u_0/q)^{M/2}
        \le\frac{\e^{1/2}}{2\sqrt2}\,\delta<\delta.
    \]
    Thus $T_\mix(P,\delta)\le t_\delta$, proving the stated bound.
\end{proof}

\section{Contraction of disagreement size}{\label{sec:contraction-disagreement}}

This section is devoted to the proof of Proposition~\ref{prop:admissible}. Let $\sigma=\frac{1+\beta_\star}{2} \in(0,1)$ be as in Lemma~\ref{lem:multi-block-boundary-comparison}, and put
\begin{equation}{\label{eq:choice-alpha-beta}}
    \beta=\min\left\{\frac1{100},\frac{\eta}{100(1+\log_{d-1}q)}\right\}, \quad \alpha=\min\left\{\frac{\beta}{100},\frac{\beta\log_{d-1}(1/\sigma)}{100} \right\}.
\end{equation}
Choose the remaining exponents as follows:
\begin{equation}\label{eq:admissible-parameters}
    \kappa=\alpha\log_{d-1}(1/\sigma), \qquad \varepsilon_0=\min\left\{ \frac{\kappa}{4},\frac14 \right\}.
\end{equation}
These are fixed positive constants, and later we will suppress all the dependencies on $\alpha, \beta, \sigma, \kappa, \varepsilon_0$. We first record a refined version of \cite[Lemma~4.3]{BGPZ26}, which controls the behavior of supercritical percolation on random regular graphs. See e.g., \cite{Pittel08, NP10} for in fact much more refined estimates.
\begin{lemma}\label{lem:product-exterior-input}
    There exist constants $K,C>0$, an integer $Q\ge0$, and a typical family $\cG_{\ref{lem:product-exterior-input}}$ with the following property. For every $G$ in this family and every $v\in V$, there is a set $H_v\subseteq S_{v,R}$, depending only on $G$, with $|H_v|\le Q$. If $Y$ is bond percolation on $G$ with parameter $\widehat p$, then:
    \begin{itemize}
        \item[(1)] For every $v\in V$ and $h\in\{r,R\}$,
        \[
            \Pb_G\Bigl( |C_{v,h}(Y)|\ge K^{-1}n,\quad |C|\le K\log n, \forall C\in\mathcal C(F_{v,h}(Y))\setminus\{C_{v,h}(Y)\} \Bigr)\ge1-n^{-30}.
        \]
        \item[(2)] Let $\xi_{Y,v}$ wire all vertices $w\in S_{v,R} \setminus H_v$ such that the component of $w$ in $F_{v,R}(Y)$ has more than $K\log n$ vertices, leaving the other vertices as singletons. For every $v\in V$ and every edge $e$ incident to $v$,
        \[
            \Pb_G\Bigl( 0\le\mu_{G[B_{v,R}]}^{\mathsf w}(e\text{ open}) -\mu_{G[B_{v,R}]}^{\xi_{Y,v}}(e\text{ open}) \le C\sigma^{(R-1)/2} \Bigr)\ge1-n^{-30}.
        \]
        In addition, for $h \in \{r,R\}$ and if $G[B_{v,h}]$ is a tree, then also
        \[
            \Pb_G\Bigl( C_{v,h}(Y)\cap S_{v,h} \text{ is }\tfrac1{4K}\text{-dense} \Bigr)\ge1-n^{-30}.
        \]
    \end{itemize}
\end{lemma}

The proof of Lemma~\ref{lem:product-exterior-input} is incorporated in \S\ref{subsec:proof-lem-5.1}. The main geometric input concerns a typical random-cluster configuration on a random regular graph. Fix a connected graph $G=(V,E)$ on $n$ vertices. Let $A \sim \mu_{G,p,q}$, and for every $v \in V$ and $A^{-} \subseteq A \subseteq A^+ \subseteq E$, we introduce the following notation:
\begin{itemize}
    \item Define $C_v=C_{v,r}(A)$, with $C_{v,r}(A)$ defined in \eqref{eq:def-fragment-cluster}.
    \item Define $R_v$ to be the largest component of $F_{v,r}(A^{-})$ in $C_v$, using the same tie-breaking rule as in \eqref{eq:def-fragment-cluster}.
    \item Recall \eqref{eq:def-local-ball}. Let $D=A^+ \setminus A^-$, denote by $\mathrm{End}(D)$ the set of vertices incident to edges in $D$, and define 
    \[
        U_v(A^{-},A,A^{+}) = \left( S_{v,r} \cap (C_v \setminus R_v) \right) \cup \left( \cup_{ \substack{ C \neq C_v: C \text{ component of } F_{v,r}(A) \\ C \cap \mathrm{End}(D) \neq \emptyset } } \big(C \cap S_{v,r} \big) \right).
    \]
\end{itemize}
\begin{prop}{\label{prop:main-geometric-property}}
    Fix $0<\alpha<\beta<1/4$ and $\varepsilon_0\in(0,1)$, and recall
    \eqref{eq:def-radii-r-R}. There exist an integer $J\ge1$ and
    a typical family $\mathcal G_{\ref{prop:main-geometric-property}}$
    such that, for every graph $G$ in this family and
    $A\sim\mu_{G,p,q}$, the following hold simultaneously with
    probability at least $1-n^{-30}$:
    \begin{itemize}
        \item[(1)] For every $v\in V$, the partition
        $\zeta_{v,r}(A)$ has at most $J$ nonsingleton blocks.
        \item[(2)] Let
        \[
            \mathfrak S(A)=\bigl\{(A^-,A^+):
            A^-\subseteq A\subseteq A^+\subseteq E,\quad
            |A^+\setminus A^-|\le\lfloor n^{1-\varepsilon_0}\rfloor
            \bigr\}.
        \]
        For every $(A^-,A^+)\in\mathfrak S(A)$,
        \[
            \sum_{v\in V}|U_v(A^-,A,A^+)|
            \le(d-1)^r(\log n+1)^{O(1)}|A^+\setminus A^-|.
        \]
    \end{itemize}
\end{prop}
We first deduce Proposition~\ref{prop:admissible} from Proposition~\ref{prop:main-geometric-property} and the local estimates of Section~\ref{sec:random-cluster-model-on-tree}. 
\begin{proof}[Proof of Proposition~\ref{prop:admissible} assuming Proposition~\ref{prop:main-geometric-property}] 
    It suffices to show that
    \begin{equation}{\label{eq:final-choice-cG}}
        \cG_{{\ref{def:typical-graph}}} \cap \cG_{\ref{lem:product-exterior-input}} \cap \cG_{\ref{prop:main-geometric-property}} \subseteq \cG_{\ref{prop:admissible}}.
    \end{equation}
    Fix $G \in \cG_{{\ref{def:typical-graph}}} \cap \cG_{\ref{lem:product-exterior-input}} \cap \cG_{\ref{prop:main-geometric-property}}$. By a simple union bound, it suffices to show that if the events in Lemma~\ref{lem:product-exterior-input} hold for $Y$ and the events in Proposition~\ref{prop:main-geometric-property} hold for $A$, then $(Y,A)$ is admissible for $T$ chosen in \eqref{eq:admissible-window}. Fix any
    \[
        Y\subseteq A^-\subseteq A\subseteq A^+\subseteq E, \qquad D=A^+\setminus A^-,\qquad k=|D|.
    \]
    We first show that
    \begin{equation}{\label{eq:bound-sum-Delta-e}}
        \sum_{e \in E} \Delta_e \le n^{-\kappa/2} |A^+\setminus A^-| +n^{-c_0}, \qquad \Delta_e:=\mu_{H_e,p,q}^{\zeta_e^+}(e\mbox{ open})-\mu_{H_e,p,q}^{\zeta_e^-}(e\mbox{ open})
    \end{equation}
    We now split the left hand side of \eqref{eq:bound-sum-Delta-e} into two parts and bound them respectively.

    \medskip
    \noindent{\bf Part 1: the summation over $e\in\mathcal B$.} Put $\zeta_e^{-}=\xi_{Y,o_e}$ and $\zeta_e^+=\mathsf w$. Clearly we have
    \[
        \zeta_e^{-} \preceq \xi_e(A^-) \preceq \xi_e(A^+) \preceq \zeta_e^+.
    \]
    Under the event in Lemma~\ref{lem:product-exterior-input}, we then have
    \[
        \Delta_e \le \mu_{H_e,p,q}^{\zeta_e^+}(e\mbox{ open})-\mu_{H_e,p,q}^{\zeta_e^{-}}(e\mbox{ open}) \le C \sigma^{(R-1)/2}.
    \]
    Since $G \in \cG_{{\ref{def:typical-graph}}}$, we have $|\mathcal B|\le dn^{3\alpha}$. Thus
    \begin{equation}\label{eq:admissible-total-cyclic-response}
        \sum_{e\in\mathcal B} \Delta_e \le Cdn^{3\alpha} \sigma^{(R-1)/2} \le n^{-2c_0}.
    \end{equation}

    \noindent{\bf Part 2: the summation over $e\notin\mathcal B$.} 
    In this case $H_e \simeq \mathbb T_{r,d}$. Write $v=o_e$, $U_e=U_v(A^-,A,A^+)$, and $W_e=C_{v,r}(Y)\cap S_e$. Under the events in Lemma~\ref{lem:product-exterior-input}, the set $W_e$ is $\rho$-dense, where $\rho=1/(4K)$. The dense set $W_e$ lies either entirely in $U_e$ or entirely in its complement. To see this, let $C$ be the component of $F_{v,r}(A)$ containing $C_{v,r}(Y)$, and let $Q$ be the component of $F_{v,r}(A^-)$ containing $C_{v,r}(Y)$. If $C=C_v$, then $Q$ either equals $R_v$ or is disjoint from it. If $C\ne C_v$, then the definition of $U_e$ includes the whole trace $C\cap S_e$ precisely when $C \cap \mathsf{End}(D) \neq \emptyset$. 

    If $W_e\cap U_e=\emptyset$, define 
    \begin{itemize}
        \item $\zeta_e^-$ to be the boundary condition obtained by
        wiring $R_v\cap S_e$ as one block, wiring $C\cap S_e$
        as a separate block for every
        $C\in\mathcal C(F_{v,r}(A))\setminus\{C_v\}$ with
        $C\cap\mathrm{End}(D)=\emptyset$, and leaving every
        remaining vertex free.
        \item $\zeta_e^+$ to be the boundary condition obtained
        from $\zeta_e^-$ by wiring $(R_v\cap S_e)\cup U_e$ as
        one block and retaining the blocks of the nongiant
        components disjoint from $\mathrm{End}(D)$.
    \end{itemize}
    A nongiant component disjoint from $\mathrm{End}(D)$ is
    unchanged in both exteriors. Every other nongiant trace
    lies in $U_e$, and $R_v$ is connected in the $A^-$-exterior.
    Moreover, every added exterior edge has both endpoints in
    $\mathrm{End}(D)$. These observations give
    \[
        \zeta_e^- \preceq \xi_e(A^-) \preceq \xi_e(A^+) \preceq \zeta_e^+.
    \]
    Moreover, the restriction of $\zeta_e^-$ and $\zeta_e^+$ on $S_e \setminus U_e$ is the same. Also, this restricted boundary condition has a block containing $W_e$ and has at most $J$ other non-singleton blocks. Lemma~\ref{lem:multi-block-boundary-comparison} gives
    \[
        \mu_{H_e}^{\zeta_e^+}(e\mbox{ open}) -\mu_{H_e}^{\zeta_e^-}(e\mbox{ open}) \le O(1)\cdot |U_e| \left( 1+J(r+1)^{O(1)} \right) \left( \frac{\sigma}{d-1} \right)^r.
    \]
    If $W_e\subseteq U_e$, instead let 
    \begin{itemize}
        \item $\zeta_e^-=\xi_{W_e}$;
        \item $\zeta_e^+=\mathsf w$.
    \end{itemize}
    Clearly we still have
    \[
        \zeta_e^- \preceq \xi_e(A^-) \preceq \xi_e(A^+) \preceq \zeta_e^+.
    \]
    Since $W_e$ is connected through the $Y$-exterior, these partitions also bracket $\xi_e(A^-)$ and $\xi_e(A^+)$. Clearly $|W_e| \ge \rho |S_e|$. Consequently, $|S_e|\le\rho^{-1}|U_e|$. Lemma~\ref{lem:dense-boundary-comparison} therefore gives
    \[
        \mu_{H_e}^{\mathrm w}(e\mbox{ open}) -\mu_{H_e}^{\xi_{W_e}}(e\mbox{ open}) \le C|S_e\setminus W_e|\left(\frac{\sigma}{d-1}\right)^r \le O(1) \cdot |U_e|\left(\frac{\sigma}{d-1}\right)^r.
    \]
    Thus, in either case we have
    \begin{equation}\label{eq:admissible-tree-response}
        0\le\Delta_e:= \mu_{H_e}^{\zeta_e^+}(e\mbox{ open}) -\mu_{H_e}^{\zeta_e^-}(e\mbox{ open}) \le O(1)\cdot |U_e| \left( 1+J(r+1)^{O(1)} \right) \left(\frac{\sigma}{d-1}\right)^r.
    \end{equation}
    Summing \eqref{eq:admissible-tree-response} over $e \notin \mathcal B$, suppose that $k\le\lfloor n^{1-\varepsilon_0}\rfloor$. Each vertex is the chosen endpoint $o_e$ for at most $d$ edges, so Proposition~\ref{prop:main-geometric-property} gives
    \[
        \sum_{e\notin\mathcal B}|U_e| \le d\sum_{v\in V}|U_v(A^-,A,A^+)| \le (d-1)^r (\log n+1)^{O(1)}|A^+ \setminus A^-|.
    \]
    This includes $k=0$, when every $U_v$ is empty. If $k>\lfloor n^{1-\varepsilon_0}\rfloor$, use $|U_e|\le|S_e|\le O(1)\cdot (d-1)^r$ instead, obtaining
    \[
        \sum_{e\notin\mathcal B}|U_e| \le Cn(d-1)^r.
    \]
    Since $\sigma^r \le \sigma^{-1} n^{-\kappa}$ and $\varepsilon_0\le\kappa/4$, the two cases imply, for all large $n$,
    \begin{equation}\label{eq:admissible-total-tree-response}
        \sum_{e\notin\mathcal B}\Delta_e \le (\log n+1)^{O(1)} n^{\varepsilon_0} \sigma^r |A^+ \setminus A^-| \le n^{-\kappa/2} |A^+ \setminus A^-|.
    \end{equation}
    Combining \eqref{eq:admissible-total-cyclic-response} and \eqref{eq:admissible-total-tree-response} yields \eqref{eq:bound-sum-Delta-e}.

    It remains to pass from these equilibrium comparisons to the local chains in Definition~\ref{def:local-chain}. Set $\delta_n=n^{-12}$. Every comparison graph has at most $M_n=O(1) \cdot n^\beta$ edges, a spanning tree of height at most $R$, and excess at most one. Every chosen comparison boundary has at most $J+1$ non-singleton blocks. Lemma~\ref{lem:spec-gap-small-excess} therefore gives, uniformly over these boundaries and over the initial configurations,
    \[
        T_{\mix}\Big( P_{H_e}^{\zeta_e^\pm},\delta_n \Big) \le O(1) \cdot M_n^2 q^{2R} \bigl(M_n+\log(1/\delta_n)\bigr) \le O(1) \cdot n^{\beta(3+2\log_{d-1} q)}.
    \]
    Choose a deterministic integer $L_n=O(1) \cdot n^{\beta(3+2\log_{d-1} q)}$ large enough for this bound, and put
    \begin{equation}\label{eq:admissible-window}
        T=\left\lceil dn\max\{2L_n,96\log n\}\right\rceil \le O(1) \cdot n^{1+\beta (3+2\log_{d-1} q)}\le n^{1+\eta/2}.
    \end{equation}
    During these $T$ full-graph update attempts, the number $N_e$ of updates retained in $H_e$ has distribution $\Bin(T,|E(H_e)|/|E(G)|)$. Its mean satisfies
    \[
        \lambda_e=\frac{T|E(H_e)|}{|E(G)|} \ge\max\{2L_n,96\log n\}.
    \]
    Chernoff's inequality gives $\Pb(N_e<L_n)\le e^{-\lambda_e/8}\le\delta_n$. Conditional on $N_e$, the retained edges are independent and uniform in $E(H_e)$, and their update uniforms remain independent. Each comparison chain is consequently within $2\delta_n$ of its stationary law after $T$ full-graph attempts. Run the lower comparison chain on $H_e$ with boundary $\zeta_e^-$ from the empty interior, and the upper one with boundary $\zeta_e^+$ from the full interior, using the same retained updates as the two actual local chains. The initial-state and boundary orders give a pathwise sandwich throughout the interval. Thus the marginal difference defining $q_e(A^-,A^+)$ is at most the marginal difference of the comparison chains. The preceding mixing bounds imply
    \[
        q_e(A^-,A^+)\le
        \begin{cases}
            \Delta_e+4\delta_n, & e\notin\mathcal B,\\
            C\sigma^{(R-1)/2}+4\delta_n, & e\in\mathcal B.
        \end{cases}
    \]
    Combining \eqref{eq:admissible-total-tree-response} and \eqref{eq:admissible-total-cyclic-response}, we conclude that
    \[
        \sum_{e\in E}q_e(A^-,A^+) \le n^{-\kappa/2} |A^+ \setminus A^-| +n^{-2c_0}+4mn^{-12} \le\frac14 |A^+ \setminus A^-|+n^{-c_0}.
    \]
    This proves that every such pair $(Y,A)$ is admissible.
\end{proof}
It remains to prove Proposition~\ref{prop:main-geometric-property}. We first state the inputs needed for our proof. The first two lemmas concern the typical behavior of random-cluster percolation on random regular graphs, whose proofs are incorporated in \S\ref{subsec:open-exterior-expansion} and \S\ref{subsec:boundary-components}, respectively.
\begin{lemma}\label{lem:reference-exterior-expansion}
    Let $r$ be as in \eqref{eq:def-radii-r-R}. There exist constants $K,\zeta>0$ and a typical family $\cG_{\ref{lem:reference-exterior-expansion}}$ with the following property. 
    %Set $a_n=K\ell_n$
    For every $G\in\cG_{\ref{lem:reference-exterior-expansion}}$ and $A\sim\mu_{G,p,q}$, the following hold simultaneously with probability at least $1-n^{-100}$:
    \begin{itemize}
        \item[(i)] For every $v\in V$, we have $|C_v|\ge\Omega(n)$ and $|C|<K\log n$ for every $C\in\mathcal C(F_{v,r}(A))\setminus\{C_v\}$.
        \item[(ii)] For an exterior component $C$, let $\partial_C W$ be the open exterior edges between $W\subseteq C$ and $C\setminus W$. For every $v\in V$ and $C\in\mathcal C(F_{v,r}(A))$,
        \begin{equation}\label{eq:reference-component-expansion}
            |W|\le O(\log n)\cdot |\partial_C W|, \qquad \emptyset\ne W\subseteq C,\quad |W|\le |C|/2.
        \end{equation}
        \item[(iii)] Every $A$-connected set $W$ satisfies
        \begin{equation}\label{eq:reference-global-expansion}
            |A \cap \partial_G W| \ge \zeta|W| \quad \text{if} \quad K\log n\le |W|\le\delta n.
        \end{equation}
        Here $\delta$ is the constant in \eqref{eq:structural-parameters}.
    \end{itemize}
\end{lemma}
\begin{lemma}\label{lem:sphere-exploration-bounds}
    Let $r$ be as in \eqref{eq:def-radii-r-R}, and let $K$ be as in Lemma~\ref{lem:reference-exterior-expansion}. There exist an integer $J_0\ge1$ and a typical family $\cG_{\ref{lem:sphere-exploration-bounds}}$ such that, for every $G$ in this family and $A\sim\mu_{G,p,q}$, with probability at least $1-n^{-100}$,
    \[
        \#\Big\{ C\in\mathcal C(F_{v,r}(A)): |C\cap S_{v,r}|\ge2, |C|<K\log n \bigr\} \le J_0 \quad \text{for every} \quad v\in V.
    \]
\end{lemma}
We also need a deterministic argument on $U_v(A^-,A,A^+)$, whose proof is postponed to \S\ref{subsec:localizing-boundary-changes}.
\begin{lemma}\label{lem:boundary-change-localization}
    Fix $\varepsilon_0\in(0,1)$ and $\alpha<\beta<1/4$, and put $R=\lfloor\beta\log_{d-1} n\rfloor$. Suppose that $G$ is $d$-regular and $\tau(G[B_{v,R}])\le1$ for every $v\in V$. Let $A\subseteq E$ satisfy the three conclusions of Lemma~\ref{lem:reference-exterior-expansion}. Then, for all sufficiently large $n$ and every $D\subseteq E$ with $|D|\le n^{1-\varepsilon_0}$, there exists a set $Z_*=Z_*(A,D)$ with $|Z_*|\le (\log n)^{O(1)} |D|$ satisfying the following condition. For every $A^-\subseteq A\subseteq A^+$ with $A^+\setminus A^-=D$, we have
    \[
        U_v(A^-,A,A^+) \subseteq S_{v,r} \cap Z_* \quad \text{for every} \quad v\in V.
    \]
\end{lemma}

We now finish the proof of Proposition~\ref{prop:main-geometric-property}.
\begin{proof}[Proof of Proposition~\ref{prop:main-geometric-property}]
    It suffices to show that $\cG_{\ref{def:typical-graph}} \cap \cG_{\ref{lem:reference-exterior-expansion}} \cap\cG_{\ref{lem:sphere-exploration-bounds}} \subseteq \cG_{\ref{prop:main-geometric-property}}$. Fix a graph $G \in \cG_{{\ref{def:typical-graph}}} \cap \cG_{\ref{lem:reference-exterior-expansion}} \cap\cG_{\ref{lem:sphere-exploration-bounds}}$. Then $G[B_{v,R}]$ has tree excess at most one. We now work on the events of Lemmas~\ref{lem:reference-exterior-expansion} and \ref{lem:sphere-exploration-bounds}, whose total failure probability is at most $2n^{-100}\le n^{-30}$. Under these events, every non-giant component $C\in \mathcal C(F_{v,r}(A))$ has fewer than $K\log n$ vertices. Lemma~\ref{lem:sphere-exploration-bounds} therefore implies that there are at most $J_0$ non-singleton blocks generated by non-giant clusters. The giant contributes at most one further block, proving Item~(1) with $J=J_0+1$.

    For Item~(2), fix $(A^-,A^+)\in\mathfrak S(A)$ and put $D=A^+\setminus A^-$. Let $Z_*=Z_*(A,D)$ be the set from Lemma~\ref{lem:boundary-change-localization}. We then have
    \begin{align*}
        \sum_{v\in V}|U_v(A^-,A,A^+)| &\le \sum_{v\in V}|S_{v,r}\cap Z_*| =\sum_{w\in Z_*}|S_{w,r}| \\
        &\le d(d-1)^{r-1}|Z_*| \le (d-1)^r(\log n)^{O(1)}|D|. \qedhere
    \end{align*}
\end{proof}

\subsection{Proof of Lemma~\ref{lem:product-exterior-input}}{\label{subsec:proof-lem-5.1}}

We first prove Item~(1). 
\begin{proof}[Proof of Lemma~\ref{lem:product-exterior-input}(1)]
    We work under the annealed setting. All estimates will apply to both $h=r$ and $h=R$. Generate the configuration model by pairing half-edges uniformly, and mark each resulting edge open independently with probability $\widehat p$, as in \cite[\S3]{NP10}. We first show that 
    \begin{equation}{\label{eq:no-intermediate-size}}
        \Pb\Big( \exists h\in\{r,R\},\ v\in V,\ C\in\mathcal C(F_{v,h}(Y)): K\log n<|C|<K^{-1}n \Big) \le n^{-200},
    \end{equation}
    and
    \begin{equation}{\label{eq:no-two-large-cluster}}
        \Pb\Big( \exists h\in\{r,R\},\ v\in V,\ C,C'\in\mathcal C(F_{v,h}(Y)): C\ne C', K^{-1}n\le |C|,|C'| \Big) \le n^{-200}.
    \end{equation}
    We first prove \eqref{eq:no-intermediate-size}. For a fixed center $v$, reveal all pairs incident to vertices at distance less than $h$. This determines $B_{v,h}$ and $S_{v,h}$ after $O(n^\beta)$ queries and leaves a uniform matching of the remaining half-edges. Whenever an exterior exploration subsequently reveals an edge with both endpoints in $S_{v,h}$, discard it, since it belongs to $E_{v,h}$. Under the event of \eqref{eq:no-intermediate-size}, there exists $w \in V \setminus B_{v,h}$ such that the component of $w$ in $F_{v,h}(Y)$ has size in $(K\log n,K^{-1}n)$. Consider the exploration process started from $w$ along its component. At each query, we match one unmatched half-edge from the current component to another half-edge, and keep it with probability $\widehat{p}$. Let $Q_j$ be the number of unmatched half-edges at vertices reached by open exterior paths after $j$ queries. Record every vertex touched by a query, including those encountered across closed edges. If the next partner belongs to a previously unrecorded vertex outside $B_{v,h}$, then $Q_{j+1}-Q_j$ equals $d-2$ with probability $\widehat p$ and $-1$ otherwise. In every other case the increment is at least $-2$. Since $\widehat{p}(d-1)>1$, choose $\lambda>0$ and $\varrho\in(0,1)$ such that
    \[
        \widehat{p}e^{-\lambda(d-2)}+(1-\widehat{p})e^\lambda<\varrho<1.
    \]
    Suppose that $N$ additional pairing queries have been made before this exploration starts. After $j$ further queries, the conditional probability that the partner has already been recorded or belongs to $B_{v,h}$ is at most
    \[
        \frac{d(Cn^\beta+2N+2j)}{dn-Cn^\beta-2N-2j-1}.
    \]
    For a sufficiently small fixed $\varepsilon_1>0$, this bound and the two possible increment estimates imply that the conditional mean of $e^{-\lambda(Q_{j+1}-Q_j)}$ is at most $\varrho$ whenever $Q_j>0$ and $j,N<\varepsilon_1n$. Let $\tau=\min\{j:Q_j=0\}$. Iterating this conditional estimate only before extinction gives
    \[
        \Eb\Big[ e^{-\lambda Q_j};\tau\ge j\Big] \le \varrho^j \qquad(1\le j\le\varepsilon_1n).
    \]
    Indeed, the induction at step $j$ uses the preceding expectation on $\{\tau>j-1\}$, which is contained in $\{\tau\ge j-1\}$; the initial factor $e^{-\lambda Q_0}$ is at most one. In particular, for every such $j$ and every $x>0$,
    \begin{equation}\label{eq:product-exterior-exploration}
        \Pb(\tau=j)\le\varrho^j,\qquad \Pb(\tau\ge j, Q_j\le xj) \le(e^{\lambda x}\varrho)^j.
    \end{equation}
    Choose $K$ large enough that $d/K<\varepsilon_1$. Observe that a complete exploration of a component with $z$ vertices takes between $z-1$ and $dz$ queries. If $K\log n<z<K^{-1}n$, its extinction time therefore lies between $\lfloor K\log n\rfloor$ and $\varepsilon_1 n$. Summing the first bound in \eqref{eq:product-exterior-exploration} over these times and over the at most $2n^2$ choices of radius, center, and root gives \eqref{eq:no-intermediate-size}.

    We then prove \eqref{eq:no-two-large-cluster}. Fix $v,h$ and put $N=\lfloor n^{2/3}\rfloor$. Under the event in \eqref{eq:no-two-large-cluster}, there exist $w_1,w_2 \in V \setminus B_{v,h}$ such that $w_1,w_2$ belong to different components in $F_{v,h}(Y)$, each of size at least $K^{-1}n$. Explore the first component until extinction or $N$ queries. The event that its size exceeds $K\log n$ but it becomes extinct before $N$ has probability at most $C\varrho^{K\log n}$. Choose $\zeta>0$ with $e^{\lambda\zeta}\varrho<1$. By the second bound in \eqref{eq:product-exterior-exploration}, the event $\{\tau\ge N,\ Q_N<\zeta N\}$ has probability at most $e^{-\Omega(N)}$. At most $O(n^\beta)$ of these half-edges are based on $S_{v,h}$. If the first component has more than $K\log n$ vertices, then outside these two failure events at least $\zeta N/2$ are based outside $B_{v,h}$. Keep those half-edges unpaired and explore the second root, unless it has already been reached by an open path in the first exploration. Stop at extinction, at an open meeting with the first component, or after $N'=\lfloor\zeta N/8\rfloor$ queries. By \eqref{eq:product-exterior-exploration}, the event that the second component has more than $K\log n$ vertices but the exploration becomes extinct before $N'$ without a meeting has probability at most $C\varrho^{K\log n}$. During these queries at least $\zeta N/4$ of the reserved half-edges remain. Each query therefore has conditional probability at least $\widehat p\zeta N/(4dn)$ of an open meeting. The resulting edge lies in the exterior because its reserved endpoint is outside $B_{v,h}$. Consequently, the probability that the two roots lie in distinct components of size at least $K^{-1}n$ is at most
    \[
        C\varrho^{K\log n}+e^{-c_1N}+e^{-c_2N^2/n}.
    \]
    A union bound over at most $2n^3$ choices of radius, center, and roots proves the two component claims with joint failure probability $O(n^{-200})$, provided $K\log(1/\varrho)>204$.
\end{proof}

We now prove Item~(2). 
\begin{proof}[Proof of Lemma~\ref{lem:product-exterior-input}(2)]
    Note that \cite[Lemma~4.3]{BGPZ26} proved that there are deterministic sets $H_{v,h}\subseteq S_{v,h}$ of at most $Q$ vertices such that the indicators of belonging to an exterior component larger than $K\log n$, indexed by $S_{v,h}\setminus H_{v,h}$, stochastically dominate independent $\Ber(1/K)$ variables. We take $H_v=H_{v,R}$ for the first inequality. The unconditioned law of $\xi_{Y,v}$ is thus $(1/K,Q)$-wired. On the component-size event in Item~(1), we have
    \[
        \xi_{Y,v} \preceq \xi_{C_{v,R}(Y)\cap S_{v,R}} \preceq \zeta_{v,R}(Y).
    \]
    Every radius-$R$ ball has excess at most one. A BFS tree has all its leaves on $S_{v,R}$, since every vertex at smaller positive depth has at least $d-2\ge1$ children. Thus these balls belong to $\mathcal U_d(R)$. By Lemma~\ref{lem:unicyclic-spatial-mixing}, applied to the unconditioned law of $\xi_{Y,v}$, the failure probability is at most $n^{-30}$ for all sufficiently large $n$, since $R\asymp\log n$. This proves the first inequality in Item~(2). In particular, this application does not condition the boundary law on the component-size event.

    Finally, if $G[B_{v,h}]$ is a tree, put $\ell=\ell_h$ as in Definition~\ref{def:dense-boundary-condition}. Each descendant subtree of height $\ell$ has $(d-1)^\ell\ge h^2$ leaves. After removing its at most $Q$ exceptional leaves, a Chernoff bound for the independent lower variables shows that at least $(d-1)^\ell/(4K)$ leaves have exterior components larger than $K\log n$, except with probability $\exp\{-c(d-1)^\ell\}$. There are at most $O(n^{1+\beta})$ such subtrees over all centers and both radii. Since $h\asymp\log n$, a union bound, including the bound $\Pb_G(\mathcal F_n)\le n^{-100}$, proves Item~(1) with probability at least $1-n^{-30}$.
\end{proof}

\subsection{Proof of Lemma~\ref{lem:reference-exterior-expansion}} \label{subsec:open-exterior-expansion}

The proof uses two estimates: connected sets of intermediate size have many open boundary edges, and deleting a small number of edges cannot produce two macroscopic components. Together with the existence of a giant component, these estimates yield the required expansion bound for every exterior component.

%By Proposition~\ref{prop:typical-graph}, there is a typical family $\mathcal H_n$ of connected $d$-regular graphs satisfying all five conditions of Definition~\ref{def:typical-graph}, both at the original radius $R$ and with $\widetilde R$ in its place. 
We take $\cG_{\ref{lem:reference-exterior-expansion}}=\mathcal G_{\ref{def:typical-graph}}$. The constants $\epsilon,\delta$ remain those in \eqref{eq:structural-parameters}. For an edge configuration $\omega$ on a graph $H$, we call a vertex set $W$ $\omega$-connected if the graph induced on $W$ by $\omega$ is connected.
\begin{lemma}\label{lem:open-boundary-expansion}
    Fix $s\in(0,1)$ such that $s(d-1-2\epsilon)>1+\epsilon$. There are constants $\zeta,c,C>0$ with the following property. Let $H$ be a graph on $n$ vertices with maximum degree at most $d$, and let $\lceil\log n\rceil\le L\le\delta n$. Suppose that every connected set $W$ in $H$ with $L\le |W|\le\delta n$ satisfies
    \[
        |E(H[W])|\le(1+\epsilon)|W|, \qquad \sum_{x\in W}(d-\mathsf{deg}_H(x))\le\epsilon|W|.
    \]
    Let $\omega\subseteq E(H)$ be a random configuration whose single-edge conditional opening probabilities are at least $s$. Then
    \[
        \Pb\Big( \exists\ \omega\text{-connected }W: L\le |W|\le\delta n, |\omega\cap\partial_H W|\le\zeta|W| \Big) \le Cn e^{-cL}.
    \]
    %The constants depend only on $d,\epsilon,s$.
\end{lemma}
\begin{proof}
    Choose $\tau>0$ so that
    \begin{equation}{\label{eq:choice-tau}}
        s(d-1-2\epsilon-\tau)>1+\epsilon.
    \end{equation}
    For sufficiently small $\lambda>0$, we have
    \[
        c_1=-\lambda(1+\epsilon)-(d-1-2\epsilon-\tau)\log(1-s+s e^{-\lambda})>0
    \]
    We then choose $0<\zeta<\min\{\tau,d/2\}$ so that $\zeta\log(ed/\zeta)<c_1/2$.

    Suppose that there exists an $\omega$-connected set $W$ with $L\le |W|=t \le \delta n$ and $|\omega \cap \partial_H W| \le \zeta|W|$. Choose an arbitrary $x \in W$ and consider the following exploration process. Start the process from $\mathtt A=\{ x \}$. Process the least labeled unprocessed edge incident to the discovered set. If the edge belongs to $\omega \cap \partial_H W$, we will mark this edge and ignore it. Otherwise, we inspect the edge and add its other endpoint to the set $\mathtt A$ if it is open. Stop when no unprocessed edge is incident to the discovered set. Since $W$ is $\omega$-connected, the discovered set is exactly $W$. Every edge incident to $W$ is processed once, so the exploration stops after $\mathtt M$ rounds, where
    \[
        dt \ge \mathtt M=\sum_{x\in W} \deg_H(x)-|E(H[W])| \ge (d-1-2\epsilon)t.
    \]
    In addition, at most $\lfloor\zeta t\rfloor$ edges are marked. This yields a sequence
    \[
        I_1,\ldots,I_{\mathtt M}, \qquad I_k \in \{ 0,1,\mathsf{mark} \}.
    \]
    For this witness, at least $\mathtt M-\zeta t$ elements of the sequence are in $\{0,1\}$, and at most $(1+\epsilon)t$ of them are $1$, since every inspected open edge lies in $H[W]$. To estimate its probability, fix the root $x$, the integers $t,\mathtt M$, and a deterministic mask $\mathcal I\subseteq\{1,\ldots,\mathtt M\}$ with $|\mathcal I|\le\lfloor\zeta t\rfloor$. Run the exploration while skipping precisely the positions in $\mathcal I$, without inspecting the skipped edges. At every other position, the next edge is determined by the preceding observations, so its conditional opening probability is at least $s$: conditioning on a partial assignment averages the full single-edge conditional probabilities. If the exploration stops early, pad the sequence with independent $\Ber(s)$ variables. Thus the $\mathtt M-|\mathcal I|$ inspected or added variables stochastically dominate independent $\Ber(s)$ variables. We enumerate deterministic masks rather than condition on the random mask associated with a witness $W$. Every witness above is produced by one of these deterministic choices.

    Since $\mathtt M-|\mathcal I|\ge(d-1-2\epsilon-\tau)t$, the exponential Markov inequality and the choice of $c_1$ give
    \[
        \Pb\bigl(\Bin(\mathtt M-|\mathcal I|,s)\le(1+\epsilon)t\bigr)\le e^{-c_1t}.
    \]
    The number of masks is at most
    \[
        \sum_{j=0}^{\lfloor\zeta t\rfloor}\binom{\mathtt M}{j}
        \le \sum_{j=0}^{\lfloor\zeta t\rfloor}\binom{dt}{j}
        \le \exp\{\zeta t\log(ed/\zeta)\}\le e^{c_1t/2},
    \]
    where the binomial entropy bound applies because $\zeta<d/2$. By a simple union bound, including the $n$ possible roots, we have
    \begin{align*}
        &\Pb\Big( \exists\ \omega\text{-connected }W: L\le |W|\le\delta n, |\omega\cap\partial_H W|\le\zeta|W| \Big) \nonumber \\
        \le\ & n \sum_{L \le t \le \delta n} \sum_{ (d-1-2\epsilon)t \le \mathtt M \le dt }
        \sum_{j=0}^{\lfloor\zeta t\rfloor} \binom{\mathtt M}{j}
        \Pb\Big( \Bin(\mathtt M-j,s) \le (1+\epsilon)t \Big) \\
        \le\ & n\sum_{L \le t \le \delta n} \sum_{ (d-1-2\epsilon)t \le \mathtt M \le dt } e^{-c_1t/2}. 
    \end{align*}
    Summing over integer $L \le t \le \delta n$ and $(d-1-2\epsilon)t \le \mathtt M \le dt$ proves the lemma, after decreasing the constant in the exponential bound.
\end{proof}
\begin{lemma}\label{lem:robust-macroscopic-uniqueness}
    Fix $G \in \cG_{\ref{def:typical-graph}}$ and recall \eqref{eq:structural-parameters}. There are constants $\gamma,c,C>0$ such that, for $A\sim\mu_{G,p,q}$, with probability at least $1-Ce^{-cn}$, every configuration $A\setminus Q$, $Q\subseteq E$, $|Q|\le\gamma n$, has at most one component with more than $\delta n$ vertices.    
\end{lemma}
\begin{proof}
    Put $h=0.01d$.
    The definition of $\epsilon$ gives $(d-1)\widehat{p}-1\ge100\epsilon$, so we may choose 
    \begin{equation}{\label{eq:choice-s-0}}
        \frac{1}{d-1} < s_0 < \widehat{p}, \qquad s_0(d-1-2\epsilon)>1+\epsilon.
    \end{equation}
    Let $\zeta,c_0,C_0$ be the constants in Lemma~\ref{lem:open-boundary-expansion} with such choice of $s=s_0$. Choose $\gamma_0>0$ and $\vartheta\in(0,1)$ such that
    \[
        0<\gamma_0\le\min\left\{ \frac{\epsilon\delta}{32}, \frac{d\delta}{400} \right\}, \qquad (1-\vartheta)\widehat{p}\ge s_0, \qquad \vartheta d/2\le\zeta\delta/16.
    \]
    First fix a deterministic $Q\subseteq E$ with $|Q|\le\gamma_0n$, and write $H=G\setminus Q$ and $\omega=A\setminus Q$. Retain each open edge of $\omega$ independently with probability $1-\vartheta$, obtaining $\omega_0$. Its conditional opening probabilities are at least $(1-\vartheta)\widehat p\ge s_0$. Moreover, conditional on the whole of $\omega_0$, the original configuration on $E(H)\setminus\omega_0$ has conditional opening probabilities at least
    \[
        \rho=\frac{\vartheta\widehat p}{1-\widehat p+\vartheta\widehat p}>0.
    \]
    Every connected set $W$ in $H$ with $\lceil\delta n/8\rceil\le |W|\le\delta n$ has $|E(H[W])|\le(1+\epsilon)|W|$ and degree deficit at most $2|Q|\le\epsilon|W|$. Lemma~\ref{lem:open-boundary-expansion}, applied to both $\omega$ and $\omega_0$, and a Chernoff bound for the number of removed open edges show that, except with probability $Ce^{-c'n}$, all three of the following hold:
    \begin{enumerate}[label=\textup{(\roman*)}]
        \item $\omega_0$ has no component whose size belongs to
        $[\lceil\delta n/8\rceil,\delta n]$;
        \item every $\omega$-connected set $W$ in this size range
        has $|\omega \cap \partial_H W|>\zeta|W|$;
        \item $|\omega\setminus\omega_0|<\zeta\delta n/4$.
    \end{enumerate}
    For the last assertion, conditional on $\omega$, the number of removed edges is dominated by $\operatorname{Bin}(dn/2,\vartheta)$. On these events, every $\omega$-component with more than $\delta n$ vertices contains an $\omega_0$-component with more than $\delta n$ vertices. Otherwise, all its $\omega_0$-components have fewer than $\lceil\delta n/8\rceil$ vertices. Contract these components and grow a connected union until its original vertex count first reaches $\delta n/2$. The resulting $\omega$-connected set $W$ has size between $\delta n/2$ and $5\delta n/8+1$, and every open edge leaving $W$ was removed by the thinning. For all sufficiently large $n$, this contradicts
    \[
        \zeta\delta n/2<|\omega \cap \partial_H W| \le|\omega\setminus\omega_0|<\zeta\delta n/4.
    \]

    Any two disjoint vertex sets of size greater than $\delta n$ are joined in $H$ by at least $c_2n$ edge-disjoint paths, each of length at most a constant $R_*$. To see this, remove at most a further $h\delta n/4$ edges and explore from the first set by breadth-first search. While the explored set avoids the second set, both it and its complement have at least $\delta n$ vertices. The expansion bound in $G$ leaves at least $h\delta n/2$ outgoing edges after the deletions. Each search step therefore adds at least $h\delta n/(2d)$ vertices until the second set is reached. We may consequently take
    \[
        R_*=\left\lceil\frac{2d}{h\delta}\right\rceil+2, \qquad c_2=\frac{h\delta}{8R_*},
    \]
    finding paths successively and removing their edges.

    Conditional on $\omega_0$, choose such paths between each pair of its components with more than $\delta n$ vertices. There are at most $\delta^{-2}$ pairs. The independent lower coupling on the remaining edges shows that the probability that a fixed pair stays disconnected in $\omega$ is at most
    \[
        (1-\rho^{R_*})^{c_2n}\le e^{-c_2\rho^{R_*}n}.
    \]
    Combining the estimates proves that, for this fixed $Q$,
    \[
        \Pb\Big(A\setminus Q\text{ has two components larger than }\delta n\Big) \le C_1e^{-c_1n}.
    \]
    Finally, choose $0<\gamma<\min\{\gamma_0,d/4\}$ sufficiently small such that
    \[
        \sum_{j=0}^{\lfloor\gamma n\rfloor}\binom{dn/2}{j}
        \le \exp\{\gamma n\log(ed/(2\gamma))+o(n)\}.
    \]
    A union bound over these sets gives the simultaneous assertion, after adjusting $c,C$.
\end{proof}

\begin{proof}[Proof of Lemma~\ref{lem:reference-exterior-expansion}]
    Fix $G \in \cG_{\ref{lem:reference-exterior-expansion}}$. By \eqref{eq:structural-parameters},
    \[
        \widehat p(d-1-2\epsilon) \ge 1+(100-2\widehat p)\epsilon>1+\epsilon.
    \]
    We may therefore apply Lemma~\ref{lem:open-boundary-expansion} with its opening-probability parameter equal to $\widehat p$. The sparsity and degree-deficit assumptions hold for $G$ and every $F_{v,r}(E)$ by Items~\textup{(4)--(5)} of Definition~\ref{def:typical-graph}; the degree deficit in $G$ is zero. Moreover, the single-edge conditional opening probabilities of $A$ and its restriction to each exterior are at least $\widehat p$, by \eqref{eq:transition-random-cluster-dynamics}. For the exterior restriction, this follows by averaging the full conditional probabilities over the edges in $E_{v,r}$. Put $a_n=\lceil K\log n\rceil$, with $K\ge1$ sufficiently large. Applying Lemma~\ref{lem:open-boundary-expansion} with $L=a_n$ and taking a union bound over these $n+1$ graphs shows that, with probability at least $1-n^{-102}$, \eqref{eq:reference-global-expansion} holds and, simultaneously for every $v\in V$,
    \begin{equation}\label{eq:open-exterior-set-expansion}
        |(A\setminus E_{v,r})\cap\partial_G W|>\zeta|W| \text{ whenever }W\text{ is connected in }F_{v,r}(A) \text{ and }a_n\le |W|\le\delta n.
    \end{equation}
    The same constant $\zeta>0$ works for the full graph and all the exteriors. In particular, none of these open graphs has a component with size in $[a_n,\delta n]$.

    We next show that $(V,A)$ has a component of size at least $a_n$ with probability at least $1-n^{-102}$. Couple $A$ with independent bond percolation $Y$ of parameter $\widehat p$ so that $Y\subseteq A$. Set
    \[
        h_n=\left\lceil \frac{\log(2a_n)}{\log((d-1)\widehat p)} \right\rceil=O(\log\log n)<r.
    \]
    Item~\textup{(2)} of Definition~\ref{def:typical-graph} gives at least $n-n^{3\alpha}$ vertices whose radius-$r$ balls are trees. A greedy selection therefore yields at least $c_d n/(d-1)^{2h_n}$ pairwise vertex-disjoint tree balls of radius $h_n$: each selected center excludes only the vertices within distance $2h_n$. In each selected ball, choose a rooted $(d-1)$-ary subtree of depth $h_n$. Let $Z_j$ be the number of vertices at depth $j$ connected to its root by $Y$. Then $Z_0=1$ and, conditional on $Z_j$, the variable $Z_{j+1}$ has law $\Bin((d-1)Z_j,\widehat p)$. The first two moment recursions give
    \[
        \Eb Z_{h_n}=((d-1)\widehat p)^{h_n}\ge2a_n, \qquad \frac{\Eb[Z_{h_n}^2]}{(\Eb Z_{h_n})^2} \le 1+\frac{1-\widehat p}{(d-1)\widehat p-1}.
    \]
    The second-moment inequality consequently gives
    $\Pb(Z_{h_n}\ge a_n)\ge c_*>0$, uniformly in $n$.
    These events are independent for the selected subtrees.
    Hence the probability that none of them contains an open
    component with at least $a_n$ vertices is at most
    \[
        \exp\left\{-\frac{c n}{(d-1)^{2h_n}}\right\}
        \le n^{-102}
    \]
    for all sufficiently large $n$. Since $Y\subseteq A$, the
    asserted existence estimate follows.

    Intersect the two events just obtained with the event in
    Lemma~\ref{lem:robust-macroscopic-uniqueness}, and let
    $\gamma>0$ be its deletion constant. Their total failure
    probability is at most $2n^{-102}+Ce^{-cn}\le n^{-100}$.
    Work henceforth on this intersection. The full-graph
    expansion estimate implies that the component just
    constructed has more than $\delta n$ vertices.
    For every $v\in V$,
    \[
        |E_{v,r}|\le C_d(d-1)^r=O(n^\alpha)=o(n).
    \]
    Deleting $E_{v,r}$ splits that component into at most
    $|E_{v,r}|+1$ components. If all of them had fewer than
    $a_n$ vertices, its original size would be at most
    $(|E_{v,r}|+1)a_n=O(n^\alpha\log n)=o(n)$, a contradiction.
    Thus every exterior has a component with at least $a_n$
    vertices, and \eqref{eq:open-exterior-set-expansion} forces
    it to have more than $\delta n$ vertices.
    Lemma~\ref{lem:robust-macroscopic-uniqueness}, applied with
    $Q=E_{v,r}$, makes this component unique. It is therefore
    $C_v$, and all other exterior components have size less
    than $a_n$, equivalently less than $K\log n$ by integrality.
    This proves Item~\textup{(i)}, while Item~\textup{(iii)}
    was already obtained from Lemma~\ref{lem:open-boundary-expansion}.

    To prove Item~\textup{(ii)}, fix $v\in V$ and an exterior
    component $C$. First suppose that $W\subseteq C$ is
    connected in $F_{v,r}(A)$ and $0<|W|\le |C|/2$. Write
    $j=|\partial_C W|$. Since $C$ is connected, $j\ge1$.
    If $|W|<a_n$, then $|W|\le a_nj$. If
    $a_n\le |W|\le\delta n$, then
    \eqref{eq:open-exterior-set-expansion} gives
    $|W|<j/\zeta$.
    It remains to consider $|W|>\delta n$. If
    $j>\gamma n/2$, then $|W|\le n<2j/\gamma$.
    Otherwise, for all sufficiently large $n$,
    \[
        |E_{v,r}\cup\partial_C W|
        \le |E_{v,r}|+j\le\gamma n.
    \]
    Deleting these edges from $A$ makes $W$ a component.
    Since $|W|>\delta n$,
    Lemma~\ref{lem:robust-macroscopic-uniqueness} implies that
    every component $W_i$ of
    $F_{v,r}(A)[C\setminus W]$ has size at most $\delta n$.
    Each $W_i$ has a nonempty boundary in $C$, consisting
    entirely of edges to $W$. The two preceding size bounds give
    \[
        |W_i|\le(a_n+\zeta^{-1})|\partial_C W_i|.
    \]
    These boundaries partition $\partial_C W$. Since
    $|W|\le |C\setminus W|$, we obtain
    \[
        |W|\le\sum_i|W_i|
        \le(a_n+\zeta^{-1})\sum_i|\partial_C W_i|
        =(a_n+\zeta^{-1})j.
    \]
    Combining all cases yields
    $|W|\le(a_n+\zeta^{-1}+2/\gamma)|\partial_C W|$.
    For an arbitrary $W\subseteq C$ with
    $0<|W|\le |C|/2$, apply this estimate to each connected
    component of $F_{v,r}(A)[W]$ and sum. Their boundaries
    are disjoint and have union $\partial_C W$, proving
    \eqref{eq:reference-component-expansion}.
\end{proof}

\subsection{Proof of Lemma~\ref{lem:sphere-exploration-bounds}} \label{subsec:boundary-components}

We will work in the annealed setting by first sampling $G \sim\cG_d(n)$ by pairing $dn$ half-edges conditioned on simplicity, and then sampling $\omega \sim \mu_{G,p,q}$. Each simple labeled $d$-regular graph has the same number $(d!)^n$ of representing pairings. Thus the joint distribution of the pairing $\pi$ and the random-cluster configuration $\omega$ is
\[
    \Pb(\pi,\omega) \propto \frac{p^{|\omega|}(1-p)^{dn/2-|\omega|} q^{k_\pi(\omega)}}{Z_\pi},
\]
where $k_\pi(\omega)$ is the number of open components and $Z_\pi$ is the FK partition function on the graph of $\pi$. Later we will slightly abuse the notation and also regard $\pi$ as a set of edges.

%We first record a switching estimate that accounts for the graph-dependent partition function in this joint law.

\begin{lemma}\label{lem:fk-weighted-switching}
    Expose pairs and their open or closed marks adaptively, choosing each next unmatched half-edge and any stopping rule using only the preceding history and independent auxiliary randomness. Let $\sigma\subseteq\pi$ be the exposed pairs and let $\mathcal F_\sigma$ be the sigma-field generated by this exposure history, including the auxiliary randomness used. Suppose that $|\pi \setminus \sigma|=N$ and $h_1$ is an $\mathcal F_\sigma$-measurable unmatched half-edge. Then for any set of unmatched half-edges $\mathcal W$ that is measurable with $\mathcal F_\sigma$, there is a constant $K_d$ such that
    \begin{equation}\label{eq:fk-weighted-partner}
        \Pb\Big( \text{the partner of } h_1 \text{ belongs to } \mathcal W \mid \mathcal F_\sigma \Big) \le\frac{2q^4|\mathcal W|}{N} \qquad(N\ge2K_d).
    \end{equation}
\end{lemma}
\begin{proof}
    From a compatible pairing containing $(h_1,h_2)$ with $h_2\in\mathcal W$, choose another unmatched half-edge $h_3$, let $h_4$ be its partner, and switch
    \[
        (h_1,h_2),(h_3,h_4) \quad\longmapsto\quad (h_1,h_3),(h_2,h_4).
    \]
    Require four distinct incident vertices and both new graph edges to be absent. At most $K_d=O_d(1)$ choices of $h_3$ are forbidden, so there are at least $N-K_d$ valid switches. Carry the first old mark to the first new pair and the second old mark to the second new pair. No observed pair or mark changes.

    Removing the two affected pairs leaves the same open graph before and after the switch. The two open component counts differ by at most two. Since the mark map is a bijection preserving the number of open edges, it follows that
    \[
        q^{-2}\le\frac{Z_{\pi'}}{Z_\pi}\le q^2, \qquad q^{-4}\le \frac{w(\pi',\omega')}{w(\pi,\omega)}\le q^4.
    \]
    An inverse switch is determined by $h_2\in\mathcal W$, because the partners of $h_1,h_2$ in the target state then determine $h_3,h_4$ and their marks. There are consequently at most $|\mathcal W|$ inverse switches per target state. Weighted double counting bounds the conditional probability in \eqref{eq:fk-weighted-partner} by $q^4|\mathcal W|/(N-K_d)$. An adaptive history, with its auxiliary random choices fixed, imposes exactly the revealed pairs and marks, so this argument applies to every permitted positive-probability history.
\end{proof}

\begin{proof}[Proof of Lemma~\ref{lem:sphere-exploration-bounds}]
    Fix $v \in V$, put $j=J_0+1$, and let $\mathcal B_v$ be the event that $S_{v,r}$ meets at least $j$ distinct components of $F_{v,r}(\omega)$, each of size less than $K\log n$, in at least two vertices each. It suffices to show that in the annealed setting 
    \begin{equation}\label{eq:small-boundary-blocks-error}
        \Pb(\mathcal B_v) \le \left[ O(1) \cdot jK(\log n) n^{-(1-2\alpha)} \right]^{j}.
    \end{equation}
    A union bound over $v \in V$ and a simple Markov's inequality then yields the desired result, by choosing $J_0$ sufficiently large such that \eqref{eq:small-boundary-blocks-error} is bounded by $n^{-300}$.

    To show \eqref{eq:small-boundary-blocks-error}, write $B=B_{v,r}$ and $S=S_{v,r}$. On $\mathcal B_v$, there are distinct components $C_1,\ldots,C_j$ of $F_{v,r}(\omega)$, each of size less than $K\log n$, and distinct vertices $u_i,w_i\in C_i\cap S$ for $1\le i\le j$. Consequently, the set $W_0=\{u_1,w_1,\ldots,u_j,w_j\}$ has size $2j$ and satisfies the event
    \[
        \mathcal E(W_0)= \Big\{ |C|<K\log n \text{ and } |C\cap W_0|\ge2 \text{ for every } C\in\mathcal C(F_{v,r}(\omega)) \text{ with } C\cap W_0\ne\emptyset \Big\}.
    \]
    Conditional on the ball-exposure history $\mathcal H_v$ described below, which determines $S$, we therefore have
    \[
        \Pb(\mathcal B_v\mid\mathcal H_v) \le \sum_{ W_0\subseteq S,|W_0|=2j } \Pb(\mathcal E(W_0)\mid\mathcal H_v).
    \]
    To bound the probability of $\mathcal E(W_0)$, first expose the $r$-neighborhood of $v$ by breadth-first search, pairing all half-edges at vertices of distance strictly less than $r$. Reveal also the open or closed mark of each exposed pair, and denote the resulting history by $\mathcal H_v$. This determines $B$ and $S$ using $O((d-1)^r)=O(n^\alpha)$ pairs. Every exposed edge lies in $B$. Since all half-edges at $B\setminus S$ have been paired, any subsequently exposed pair is internal to $B$ if and only if both endpoints lie in $S$. Fix a realization of $\mathcal H_v$. Start from the least labeled vertex of $W_0$ not contained in a previously completed component, and explore its open component in $F_{v,r}(\omega)$. Process incident half-edges at the vertices reached from the root by open exterior edges, using a fixed order. Whenever a new pair is queried, reveal its partner and its mark. Ignore pairs with both endpoints in $S$ when determining exterior connectivity. Stop the entire procedure as soon as the current component reaches $K\log n$ vertices. Otherwise, complete that component and repeat until every vertex of $W_0$ lies in a completed component. At most $2j$ components are started, and hence the procedure makes at most $M=2dj\lceil K\log n\rceil$ new pairing queries, whether or not $\mathcal E(W_0)$ occurs.

    Let $I_t$ indicate that the partner revealed in the $t$-th new query is incident to a vertex of $W_0$, irrespective of the mark or whether the pair is internal. Put $I_t=0$ if the procedure has already terminated. On $\mathcal E(W_0)$, every explored component has fewer than $K\log n$ vertices, so the procedure never stops early. Each completed component contains at least two vertices of $W_0$, so at most $j$ vertices of $W_0$ are used as roots. Each of the other at least $j$ vertices must occur as the partner in a new pairing query: before being reached by an open exterior edge, it cannot be processed as a query source unless it is chosen as a root. All edges revealed in the ball exposure are internal, and any open exterior edge revealed later is traversed when it is revealed. Distinct such vertices require distinct queries. We have therefore proved the deterministic implication
    \[
        \mathcal E(W_0)\subseteq
        \left\{\sum_{t=1}^{M}I_t\ge j\right\}.
    \]

    We now bound the probability of this event. Throughout the procedure, only $O(n^\alpha)+M=o(n)$ pairs have been exposed, so the number of unexposed pairs remains a positive multiple of $n$. At any query, the set of unmatched half-edges incident to $W_0$ is measurable from the preceding history and has size at most $2dj$. Lemma~\ref{lem:fk-weighted-switching} therefore gives
    \[
        \Pb(I_t=1\mid\text{preceding history})\le\frac{Cj}{n}.
    \]
    For any prescribed $j$ query positions, iterated conditioning bounds the probability of a hit at every one of them by $(Cj/n)^j$. Taking a union bound over these positions yields
    \[
        \Pb(\mathcal E(W_0)\mid\mathcal H_v)
        \le \binom{M}{j}\left(\frac{Cj}{n}\right)^j
        \le \left(\frac{C'jK\log n}{n}\right)^j.
    \]
    Since $|S|\le C_dn^\alpha$, there are at most $(C_dn^\alpha)^{2j}$ choices of $W_0$. A union bound over these choices, followed by averaging over $\mathcal H_v$, gives \eqref{eq:small-boundary-blocks-error}. Choosing $J_0$ large enough makes the annealed probability that any center violates the bound at most $n^{-299}$. Markov's inequality then shows that, outside a set of graphs of probability at most $n^{-199}$, the conditional failure probability is at most $n^{-100}$.
\end{proof}

\subsection{Proof of Lemma~\ref{lem:boundary-change-localization}} \label{subsec:localizing-boundary-changes}

Put $h=R-r\asymp\log n$. We first record two elementary observations. First, fix any nonempty connected vertex set $T$ in $F_{v,r}(E)$. Consider the components of $G[B_{v,R}]\setminus E_{v,r}$ that meet $S_{v,r}$. The spanning-tree argument in the proof of Lemma~\ref{lem:small-ball-removal-degree} gives
\[
    \sum_{C:\text{such components}} \Big( |C\cap S_{v,r}|-1 \Bigr) \le\tau(G[B_{v,R}]) \le 1.
\]
After discarding at most one vertex of $T\cap S_{v,r}$, the remaining sphere vertices lie in distinct components of this annular graph. An exterior path between two of them must leave $B_{v,R}$ and therefore has length at least $2(h+1)$. Their radius-$h$ balls within the connected exterior graph on $T$ are consequently disjoint. Each contains at least $\min\{|T|,h+1\}$ vertices, so
\begin{equation}\label{eq:exterior-sphere-spacing}
    |T\cap S_{v,r}| \le 1+\frac{|T|}{\min\{|T|,h+1\}} \le 2+\frac{|T|}{h+1}.
\end{equation}
Next, in any graph $H$ of maximum degree at most $d$, for a fixed vertex $x$ and a fixed integer $j\ge0$,
\begin{align}
    & \#\Big\{ T:x \in T, T\text{ is connected in }H, |T|<K\log n, |\partial_H T|\le j \Big\} \nonumber \\
    \le\ & \sum_{i=0}^j\binom{dK\log n}{i} =O_{j}((\log n)^j).  \label{eq:small-connected-set-count}
\end{align}
To see this, explore from $x$, processing incident edges in a fixed deterministic order. To discover exactly $T$, skip its boundary edges and traverse every other processed edge. At most $d|T|<dK\log n$ edges are processed, and at most $j$ positions are skipped. The root and the skipped positions determine the exploration and hence $T$, proving the bound.

\begin{proof}[Proof of Lemma~\ref{lem:boundary-change-localization}]
    Fix $D$ and write $k=|D|$. If $k=0$, take $Z_*=\emptyset$. Assume henceforth that $k\ge1$, and let $X=\mathsf{End}(D)$. Starting from $Z=X$, repeatedly add any connected set $T$ in the remaining graph $A[V\setminus Z]$ such that
    \[
        K\log n\le |T|\le\delta n/2, \qquad |\partial_{A[V\setminus Z]}T| \le\frac{\zeta}{4}|T|.
    \]
    (Use a fixed deterministic rule to choose among eligible sets, and stop when no such set remains.)
    Each addition decreases the open boundary of $Z$.
    More precisely, writing $c=|\partial_{A[V\setminus Z]}T|$
    and using \eqref{eq:reference-global-expansion},
    \[
        |A\cap\partial_G(Z\cup T)|-|A\cap\partial_G Z| =2c-|A\cap\partial_G T| \le-\frac{\zeta}{2}|T|.
    \]
    The added sets are disjoint. Summing this inequality therefore gives
    \begin{equation}\label{eq:pruned-set-size}
        |Z| \le |X|+\frac{2}{\zeta}|A\cap\partial_G X| \le\left(1+\frac{2d}{\zeta}\right)|X| =O(k).
    \end{equation}
    Let $H$ be the open graph induced by $A$ on $V\setminus Z$. By termination,
    \begin{equation}\label{eq:pruned-graph-expansion}
        |\partial_H T|>\frac{\zeta}{4}|T|
        \text{ for every connected }T\subseteq H
        \text{ with } K\log n\le |T|\le\delta n/2.
    \end{equation}

    Now fix $A^-\subseteq A\subseteq A^+$ with $A^+\setminus A^-=D$, a vertex $v \in V$, and a vertex $w\in U_v(A^-,A,A^+)\setminus Z$. If $w\in C_v\setminus R_v$, let $W$ be its component in $F_{v,r}(A^-)$.
    Since $R_v$ is the largest retained component of $C_v$,
    we have $|W|\le |C_v|/2$.
    Every open exterior edge leaving $W$ belongs to
    $A\setminus A^-\subseteq D$, and at least one such edge
    exists because $W$ is a proper subset of the connected
    component $C_v$.
    Thus $W\cap X\ne\emptyset$, and
    \eqref{eq:reference-component-expansion} gives
    $|W|\le O(k\log n)$.
    Otherwise, let $W\ne C_v$ be the component of
    $F_{v,r}(A)$ containing $w$.
    By the definition of $U_v$ and Item~(i) of
    Lemma~\ref{lem:reference-exterior-expansion},
    we again have $W\cap X\ne\emptyset$ and $|W|<K\log n$.
    In both cases, $W$ is connected in $F_{v,r}(A)$ and
    \[
        W\cap X\ne\emptyset,\qquad |W|\le O(k\log n)=o(n).
    \]
    Let $T$ be the component containing $w$ in
    $H\setminus E_{v,r}$. All endpoints of $D$ belong to $Z$,
    so $A$ and $A^-$ agree on edges whose endpoints lie outside
    $Z$. Every exterior $A$-edge leaving $W$ belongs to
    $D$, so $T\subseteq W\setminus Z$.
    Moreover, every edge of $\partial_H T$ belongs to $E_{v,r}$. Every vertex of $B_{v,r}\setminus S_{v,r}$ is isolated in $H\setminus E_{v,r}$. Since $w\in S_{v,r}$, its component $T$ therefore avoids $B_{v,r}\setminus S_{v,r}$. Consequently, the endpoint in $T$ of every edge of $\partial_H T$ belongs to $T\cap S_{v,r}$.
    Equation~\eqref{eq:exterior-sphere-spacing}
    therefore gives
    \[
        |\partial_H T|
        \le d|T\cap S_{v,r}|
        \le d\left(2+\frac{|T|}{h+1}\right).
    \]
    Note that $|T|\le |W|\le\delta n/2$.
    In addition, if $|T|\ge K\log n$, the preceding bound contradicts
    \eqref{eq:pruned-graph-expansion}, since
    \[
        \frac{|\partial_H T|}{|T|}
        \le\frac{2d}{K\log n}+\frac{d}{h+1}
        =o(1)<\frac{\zeta}{4}.
    \]
    Hence $|T|<K\log n$. Since $K\log n/(h+1)=O(1)$, there is a fixed
    integer $J_1$ such that
    \[
        |\partial_H T|
        \le d\left(2+\frac{K\log n}{h+1}\right)\le J_1
    \]
    for all sufficiently large $n$.
    Finally, $T$ contains a vertex adjacent to $Z$ through an
    open edge. Indeed, follow an exterior path in $W$ from
    $w$ to $W\cap X$; before its first visit to $Z$, the path
    remains in $T$.

    Define $Z_*$ by adjoining to $Z$ every connected set $T$
    in $H$ with $|T|<K\log n$ and $|\partial_H T|\le J_1$ that
    contains a vertex adjacent to $Z$ in $A$.
    There are at most $d|Z|$ possible choices of this vertex.
    Equations~\eqref{eq:small-connected-set-count} and
    \eqref{eq:pruned-set-size} imply
    \[
        |Z_*|
        \le |Z|+d|Z| \cdot (K\log n) \cdot O((K\log n)^{J_1})
        \le k (\log n)^{O(1)}.
    \]
    The preceding argument shows that it contains every vertex
    of $U_v(A^-,A,A^+)$, for every center $v$ and every admissible
    pair with difference set $D$. Since $U_v\subseteq S_{v,r}$,
    this proves the lemma.
\end{proof}

\subsection*{Acknowledgments}

We thank Hang Du for stimulating discussions at an early stage of the project. This work is partially supported by the NSFC Major Program (Project No. 12595294), the NSFC Key Program (Project No. 12231002), and the New Cornerstone Investigator Program (Project No. NCI202501).

\subsection*{Statement on AI use}

After the work \cite{LS26}, the author realized that it might be possible to prove Theorem~\ref{thm:main} using a similar bootstrap argument for the disagreement size. Building on this strategy, the main results of this paper were obtained using GPT-6 and the Pharos system. Pharos is an AI-powered mathematical research system built upon and substantially extending Rethlas~\cite{JGJ+26} and Danus \cite{LGS+26}, which is publicly available on GitHub and was first introduced in \cite{DGG+26}. We thank Bin Dong, Guoxiong Gao, Jiedong Jiang, Shurui Liu, Zeming Sun, and Bin Wu for developing and refining Pharos.

During the subsequent verification and revision, the author substantially simplified the argument. One important simplification was in the proof of Proposition~\ref{prop:main-geometric-property}, which controls the cumulative effect of disagreements on the local boundary conditions. Proposition~\ref{prop:main-geometric-property} was originally proved using delicate control of the spectral properties of killed random walks on random regular graphs. The author realized that this proposition can be proved using more straightforward structural estimates for the supercritical random-cluster model on random regular graphs, which was later formalized in Lemmas~\ref{lem:reference-exterior-expansion}--\ref{lem:boundary-change-localization}. All AI proofs in this paper were digested and rewritten by the author, and the author takes full responsibility for the results.

\bibliographystyle{alpha}

\end{document}